\documentclass[10pt]{article}

\usepackage{amsmath}
\usepackage{mathtools}
\usepackage{amsthm}
\usepackage{authblk}
\usepackage{color}
\usepackage{amssymb}
\usepackage{amsfonts}
\usepackage[graphicx]{realboxes}
\usepackage{verbatim}
\usepackage{mathrsfs}
\usepackage{subfigure}
\usepackage{fancyhdr}
\usepackage{latexsym}
\usepackage{dsfont}
\usepackage{fancybox}
\usepackage{bbm}
\usepackage{pdflscape}

\usepackage{booktabs}
\usepackage{multirow}
\usepackage{adjustbox}

\usepackage{framed, color}
\definecolor{shadecolor}{rgb}{1, 0.8, 0.3}
\usepackage{bm}

\def\bR{\mathbb{R}}
\def\cco{{\bf C}^0}

\def\ccol{{\bf C}^{0,\ell}}
\def\ccoL{{\bf C}^{0,L}}
\def\cco1{{\bf C}^{0,1}}
\def\bP{\mathbb{P}}

\def\o{{\bf {Ope}}}
\def\o{\mathbb{O}}
\def\jj{{\mathbb{J}}}

\newtheorem{lemma}{Lemma}

\newtheorem{proposition}{Proposition}

\newtheorem{cor}{Corollary}
\definecolor{lightgreen}{rgb}{0.8,1,0.8}
\definecolor{darkgreen}{rgb}{0,0.7,0}
\definecolor{lightred}{rgb}{1,0.8,0.8}
\definecolor{darkred}{rgb}{0.9,0,0}
\definecolor{lightblue}{rgb}{0.8,0.8,1}
\definecolor{darkblue}{rgb}{0,0,0.6}
\definecolor{lightyellow}{rgb}{1,1,0.7}

\newcommand{\colb}[1]{{\bf\color{darkblue} \1}}
\newcommand{\colg}[1]{{\bf\color{darkgreen} \1}}
\newcommand{\colr}[1]{{\bf\color{darkred} \1}}

\title{
State dependent delay maps: numerical algorithms 
and dynamics of projections. 
}

\author{J.D. Mireles James \thanks{J.M.J partially 
supported by NSF grant DMS - 1813501
Email: {\tt jmirelesjames@fau.edu}}}

\author{Francis Motta \thanks{ 
Email: {\tt fmotta@fau.edu}}}

\author{Vincent Naudot \thanks{
Email: {\tt vnaudot@fau.edu}}}

\affil{Florida Atlantic University, Department of Mathematical Sciences,
777 Glades Road, 33431 Boca Raton, FL, USA}
\date{\today}

\begin{document}
\maketitle

\begin{abstract}
This work concerns the dynamics of a certain class of
delay differential equations (DDEs) which we refer to as state dependent
delay maps.  These maps are generated by delay differential equations
where the derivative of the current state depends only on delayed variables,
and not on the un-delayed state.  However, we allow that the delay
is itself a function of the state variable.  
A delay map with constant delays can be rewritten explicitly 
as a discrete time dynamical system on an appropriate function space, and a 
delay map with small state dependent terms can be viewed as a ``non-autonomous'' 
perturbation. We develop a fixed point formulation for the Cauchy 
problem of such perturbations, 
and under appropriate assumptions obtain the existence of forward 
iterates of the map.  

The proof is constructive and leads to numerical procedures
which we implement for illustrative examples, including the cubic Ikeda and 
Mackey-Glass systems with constant and state-dependent delays. 
After proving a local convergence result for the method, we
study more qualitative/global convergence issues using data analytic tools for time 
series analysis (dimension and topological measures derived 
from persistent homology).  Using these tools we quantify the 
 convergence of the dynamics in the finite dimensional projections
 to the dynamics of the infinite dimensional system.  
 \end{abstract}

\section{Introduction} \label{sec:intro}
Existence, uniqueness, and regularity theory for state dependent delay 
differential equations (SDDDEs) is notoriously delicate, in part because 
the natural setting for the Cauchy problem is an infinite dimensional 
Banach manifold 
\cite{MR4232666,MR2457636,MR2084753,MR2019242,MR3642771,MR4261206,MR4396234}.
The resulting difficulties are amplified when considering the 
qualitative dynamics and, in spite of much progress,
there are many important open questions.  
We refer to the works of
\cite{MR3773293,MR4222179,MR3537351,MR3547461}
for more thorough discussion of qualitative topics like invariant 
manifold theory in this context.

One alternative to studying the orbit structure generated by the 
Cauchy problem on the full phase space (whatever that might be
in the case of SDDDEs), is to surgically focus instead on isolated,
geometrically meaningful special solutions or collections of solutions.
Examples include equilibria, periodic solutions, 
quasi-periodic solutions, stable/unstable manifolds, 
etcetera, and we exploit the fact that such landmark solutions 
may enjoy better regularity properties than typical orbits.
More precisely, these landmarks can be reformulated as solutions of 
invariance equations with nicer functional analytic properties than 
the original SDDDE Cauchy problem.  
Moreover, it is often possible to develop a-posteriori analytical methods 
for the functional equations describing a landmark.  
Recent successful examples of this 
approach include the work of
\cite{gimeno2021persistence,MR4308692,Per,MR4287353}
on periodic orbits and their whiskers, and the work of
\cite{MR4112213,MR3501842,MR3736145} on both KAM and hyperbolic invariant 
tori.

In the present work we consider an a-posteriori approach to the Cauchy problem for 
SDDDEs which arise as perturbations of delay differential equations (DDEs) with constant delays.  
We remark that the quantitative and qualitative theories for constant delays are well developed,
and refer for example to the  book of \cite{MR1243878}.
The idea of our approach is to treat long enough solutions 
of the constant delay problem
as approximate solutions of the state dependent perturbation.  
In the present work we further restrict our 
attention to a simple class of constant delay systems where the 
right hand side of the equation depends only on the past.  
This simplifies the unperturbed problem dramatically, in the sense
that we can express it as an explicit discrete time semi-dynamical system
on an appropriate function space.  
Our method is iterative, 
and we show that it converges under appropriate conditions
to a solution of the SDDDE.

In the second part of the paper, we are interested in more global 
questions about the robustness of attractors with respect to perturbation.
In this later portion of the paper we do not prove any theorems, 
but apply and compare several 
complexity indicators from topological data and time series analysis.  
These data analytic tools are used to measure the complexity of the attractors in the 
truncated system as a function of the projection dimension.  We observe 
rapid convergence of the iterative scheme to the ``true dynamics'' of the system, 
where ground truth is measured by taking very high projection dimension.  
The results provide quantitative measures of the size of the projection
dimension needed to accurately describe the SDDDE.

\section{Formal problem statment} \label{sec:setup}

Let $F \colon \mathbb{R} \to \mathbb{R}$ be a smooth map and $\tau > 0$.
In practice, $F$ only needs to be defined on some open subset of $\mathbb{R}$.
Consider the constant delay differential equation 
\begin{equation}\label{1b}
y'(t) = F(y(t - \tau)).
\end{equation}
Note that the right hand side depends only on the history of the system
(in a more general problem -not considered here - $F$ would depend also on $y(t)$).

The Cauchy problem for Equation \eqref{1b} is as follows.  
Suppose that a continuous function $y_0 \colon [-\tau, 0] \to \mathbb{R}$ is given.
Find a $T>0$ and a function $y \colon [-\tau, T) \to \mathbb{R}$ having 
\[
y(t) = y_0(t) \quad \quad \quad t \in [-\tau, 0], 
\]
which solves Equation \eqref{1b} on $(0, T)$.
Note that for Equation \eqref{1b} the existence of such a function for $T = \infty$ is 
trivial.  One simply defines $y_1 \colon [0, \tau] \to \mathbb{R}$ by 
\[
y_1(t) = y_0(0) + \int_0^t F(y_0(s-\tau))\, ds,
\]
and iteratively 
\[
y_{n+1}(t) = y_{n}(n\tau) + \int_{n\tau}^t F(y_n(s-\tau)) \, ds.
\]
For any $N \in \mathbb{N}$, 
the function $y \colon [-\tau, N\tau] \to \mathbb{R}$ defined by 
\[
y(t) = \begin{cases}
y_0(t) & -\tau \leq t < 0 \\
y_1(t) & 0 \leq t < \tau \\
& \vdots \\
y_N(t) & (N-1)\tau \leq t < N \tau 
\end{cases},
\]
solves the Cauchy problem.  In fact, nothing prevents us from taking $T \to \infty$.  
The solution may be unbounded as $t \to \infty$, but it exists on $[-\tau, \infty)$.
Note that, since $y_0$ is continuous, $y_1$ is $C^1$ and $y_n$ is $C^n$.  
This argument is known as the method of steps, and in the context of the simple 
DDE of Equation \eqref{1b} it can be thought of as a discrete time infinite 
dimensional (semi) dynamical system on \(C([-\tau, 0], \mathbb{R})\) 
(i.e., the space of real valued continuous function defined on the interval \([-\tau,0]\)).

The situation becomes much more complicated when we allow the delay to 
depend on the state of the system.  So, let $\delta \colon \mathbb{R} \to \mathbb{R}$
be a smooth function and consider the delay differential equation 
\begin{equation}\label{1a}
y'(t) = F\left( 
x(t - \delta(x(t)))
\right).
\end{equation}
Note that when $\delta \equiv \tau$ is constant we recover the earlier case.
In this work we restrict our attention to the perturbative situation where
\[
\delta(u) = \tau + \varepsilon h(u). 
\]
Indeed, we consider the case when $h(u) = u$ is the identity mapping.  

Analyzing the resulting Cauchy problem is more subtle in this case.  Given a 
 continuous (or smooth) function $y_0 \colon [-\tau, 0] \to \mathbb{R}$ can we 
 find a solution of Equation \eqref{1a} on the interval $[0, \tau]$?
 Rerunning the method of steps argument shows that the answer is delicate.
 We may in fact need that $y_0$ is defined on \([-\tau - {\hat \delta}, 0]\), where 
 \({\hat \delta}\) depends on the supremum norm of \(y_0\).  This fact --
 that the domain of the history function depends on the 
 size of the history function -- leads to the  
  phase space for Equation \eqref{1a} being a Banach manifold.  
 
 In Section \ref{sec:stateDependentCase} we study this problem in more detail, 
  working near a solution of the 
 constant delay problems which has as much history as we need.  That is,
 we only try to follow orbits of the constant delay system which exist for 
 sufficiently long time into the state dependent perturbation.  The intuition behind this
 is that we want to study the perturbed problem near an attractor of the 
 constant delay system, and solutions on an attractor have infinite prehistory.

\bigskip

Effectiveness of the iterative numerical scheme is examined in the 
case of two explicit example problems.
These are:

\begin{itemize}
\item \textbf{The Cubic Ikeda family:} 
\[
\dot y =  y(t-\delta(y(t))) - y^3(t-\delta(y(t))),\ \ {\mathrm {i.e.}}\ \ F(u) = u -u^3.
\]
When $\varepsilon = 0$, this system was studied extensively in \cite{MR2570698},
where simulation results are presented which suggest the existence of a chaotic
attractor.
\item \textbf{The Mackey-Glass Family:} 
consider
\[
\dot y =  -a y + \beta\frac{y(t-\delta(y(t)))}{1+y^n(t-\delta(y(t)))},
\ \ {\mathrm {i.e.}}\ \ F(u) = \beta\frac{u}{1+u^n}
\]
where $\beta, \ n>0$.  When $\varepsilon = 0$ the system was originally studies in 
 \cite{1977Sci...197..287M}, and was one of the first infinite dimensional dynamical systems 
 conjectures to have a strange attractor/chaotic dynamics.  
\end{itemize}

Note that Mackey-Glass is not precisely of the form given in Equation \eqref{1b}.  Rather 
it has a linear ``friction'' term $-a y(t)$.  Nevertheless, it is well known that such a linear term 
can be removed, transforming the system into the appropriate form.  This is 
reviewed in Section \ref{sec:friction}.

\section{The state dependent case} \label{sec:stateDependentCase}
We rescale the delay in equation (\ref{1b})  in the following manner.
Write
$x(t) =y(\tau t)$, from (\ref{1b}) we  get
\begin{equation}\label{2}
\dot x =\tau F[x(t-\frac{1}{\tau}\delta(x(t))) ].
\end{equation}
Since
$
\delta(u)/\tau = 1 -  \varepsilon u
$,  
 (\ref{2}) finally writes
\begin{equation}\label{1}
\dot x =  \tau F[x(t-1+{\varepsilon} x(t)) ].
\end{equation}
We develop an algorithm to compute solutions of 
(\ref{1}) by perturbing  away from 
the following phase independent delay equation,
\begin{equation}\label{3}
\dot x = \tau F[ x(t-1)], \ \ \ \ \ \ \ \ t\geq 0 
\end{equation}
when $\varepsilon = 0$.


Recall that a solution $x=x(t)$ of (\ref{1}) needs to satisfy
\begin{equation}\label{1bis}
\left\{
\begin{array}{rcl}
\displaystyle\frac{dx}{dt}& = &  \tau F\biggl(x(t-1+\varepsilon x(t))\biggr),  \ \ \ t \geq 0 \\
\\ & \ & \\
x(t) & = &  x_0(t), \  \ \ \ \ \ \ \ t \in [-3/2,0]
\end{array}
\right.
\end{equation}
Our interest in this article is for bounded orbits. Therefore, we shall assume that 
there exists a contant \(M>1\) such that
the 
solution of (\ref{1bis}) satisfies
\begin{equation}\label{aboutMM}
|x(t)| \leq M,  \ \ \forall t \geq -3/2.
\end{equation}

Let \({\mathrm {I}} \subset {\bR}\) be a compact interval and \(\ell>0\).
We say that a continuous function 
\({\bf f} \in C({\mathrm {I}},{\bR})\)
belongs to
\(
\ccol({\mathrm {I}})
\) 
if \({\bf f}\) 
is diferentiable on \({\mathrm {int}}({\mathrm {I}})\) the interior of \({\mathrm {I}}\)
with
\[
  |{\bf f}^\prime(t)| \leq \ell,  \ \ \forall \ t\in {\mathrm {int}}({\mathrm {I}}). 
\]
Our 
goal is to construct a solution $x$ of (\ref{1bis}) defined for all positive real number such that
$x$ coincides with $x_0$ on 
\([-3/2,0]\)
where $x_0$ 
satisfies (\ref{aboutMM})
and belongs to  $\ccol([-3/2,0]) $ 
for some \(\ell\) large enough 
(see below).
Our approach is to extend the definition of $x$ on a larger domain, by showing that this 
extension 
satisfies a fixed point argument. 

\subsection{A contraction }
Let \(M>1\) be a (large) real number. 
We define
\begin{equation}\label{bnd3}
    K_0=\sup_{|\xi| \leq M} |F(\xi)|, \ \ K_1=\sup_{|\xi| \leq M} |F'(\xi)|, \ \ 
    K_{2}=M+\tau K_{0}/2.
    \end{equation}
Let \(\ell \geq \tau K_0\)
and $z \in \ccol([-3/2, 0])$  such that for all \(-3/2\leq t\leq 0\),
\begin{equation*}
|z(t)| \leq M.
\end{equation*}
We define
\begin{equation}\label{k22}
\ccol_0([0,1/2])
  = \{ {\bf f} \in  \ccol([0,1/2]), 
  \ | \ {\bf f}(0)=z(0) , \ \  \sup_{0\leq t\leq 1/2} |{\bf f}(t)|
  \leq K_{2} \}.
\end{equation}
This set
\(\ccol_0([0,1/2])\)
is equippied with the following distance
\begin{equation}\label{distt}
    {\bf d}({\bf f}, {\bf g})=
  \sup_{t \in (0,1/2)} |{\bf f}^\prime(t)-{\bf g}^\prime(t)|.
\end{equation}
and the following operator 
\[
\o_z: \ccol_0([0,1/2])\rightarrow \ccol_0([0, 1/2]), \ \ w \mapsto \o_z(w) 
\]
where
\begin{equation}\label{oo}
\o_z(w)(t) = z(0) + \tau \int_{0}^{t} F\biggl(z(s-1+\varepsilon w(s))\biggr)ds.
\end{equation}
For convenience, we also use the following notation
\[
|{\bf g}|_{0}=\sup_{0\leq t \leq 1/2}|{\bf g}(t)|, \ \ |{\bf g}|_{1}=\sup_{0< t < 1/2}|{\bf g}^{\prime}(t)|.
\]

\noindent
We state the following proposition.
\begin{proposition}\label{p1}
{\sl 
There exists $\varepsilon_0>0$ such that for all
$|\varepsilon|\leq \varepsilon_0$,
the operator $\o_z$ is a contraction. More precisely
for all $w_1, \ w_2 \in  \ccol_0([0,1/2])$, 
\[
  {\bf d}
  \biggl(
  \o_z(w_1),
  \o_z(w_2)
  \biggr)
\leq \frac{1}{2} 
  {\bf d}(w_1, w_2)
  .
\]
}
\end{proposition}
\bigskip
\noindent
Observe that $\ccol_0([0,1/2])$ equipped with the distance defined in 
(\ref{distt})
is a complete space.
As a consequence,
for each \(z \in \ccol([-3/2,0])\), \(\o_z\) admits a unique fixed point 
\({ \Psi}_\varepsilon(z)\in \ccol_0([0,1/2])\) i.e., 
satisfies
\[
{\Psi}_\varepsilon(z)(t)=z(0) + \tau \int_{0}^{t} F\biggl(z(s-1+\varepsilon { \Psi}_\varepsilon(z)(s))\biggr)ds.
\]
The map
\[
 \Psi_\varepsilon:    
\{ {z} \in \ccol([-3/2,0]) \ | \ \sup_{-3/2\leq t\leq 0}|{z}(t)| \leq M\} 
\rightarrow \ccol_0([0,1/2]),  \ z \mapsto {\Psi}_\varepsilon(z)
\]
is called the Half-Stroboscopic map.
Furthermore, (as mentionned in the introduction), if 
\(z\) is \({\bf C}^r\), 
\({\Psi}_\varepsilon(z)\) is \({\bf C}^{r+1}\) on \((0,1/2)\) and we have
\[
\frac{d}{dt}{\Psi}_\varepsilon(z)(t) = \tau F\biggl(x(t-1+\varepsilon {\Psi}_\varepsilon(z)(t))\biggr).
\]
We now apply Proposition 1 to our initial condition \(x_0\).
The function
\[
x_{1/2} :[-2, 0 ] \rightarrow  {\bR},  \ t \mapsto 
x_{1/2}(t)
\]
where
\begin{eqnarray*}
x_{1/2}(t) & = &  x_0(t+1/2)  \  \ \ \ \ \ \ 
\ \ \ {\mathrm {if}} \     -2\leq t\leq -1/2 ,
\\ 
& \ & \\
\ & = &  {\Psi}_\varepsilon(x_0)(t+1/2)  
 \ \ \  \ {\mathrm {if}} \    - 1/2\leq t\leq 0
\end{eqnarray*}
is \(\ell\)-Lipschitz. 
It is easy to verify that the function
\begin{eqnarray*}
x(t) & = &  x_{1/2}(t-1/2),
\end{eqnarray*}
is well
defined on the interval \([-3/2,1/2]\) 
 and is a solution of (\ref{1bis}).
 We construct the sequence of functions
 \[
    (x_{m/2}) , \ m\geq 1
    \]
    (and deduce thereafter 
    the solution of 
(\ref{1bis}) for all real number)
by induction on \(m\)  in the following way. 
At each step of the induction, for \(m\geq 0\), we construct \(x_{m/2}\) from \(x_{(m-1)/2}\)
assuming the latter to satisfy
(\ref{aboutMM}), more precisely 
\begin{equation}\label{aboutMMM}
|x_{(m-1)/2}(t)| \leq M \ \  \ \ \ {\mathrm{for}}\ {\mathrm {all}} \ \ -\frac{(2+m)}{2}\leq t<0.
\end{equation}
We 
then define 
\[
x_{m/2}:[-3/2-m/2,0]  \rightarrow  {\bR},  \ t \mapsto
x_{m/2}(t)
\]
where
\begin{eqnarray*}
x_{m/2}(t) & = &  x_{(m-1)/2}(t+1/2)  \  \ \ \ \ \ \ 
\ \ \ {\mathrm {if}} \     -3/2-m/2\leq t\leq -1/2 ,
\\ 
& \ & \\
\ & = &  {\Psi}_\varepsilon(x_{(m-1)/2})(t+1/2)   
 \ \ \  \ {\mathrm {if}} \    - 1/2\leq t\leq 0,
\end{eqnarray*}
and 
again, it is easy to verify that
the function
\begin{eqnarray*}
x(t) & = &  x_{m/2}(t-m/2),
\end{eqnarray*}
is well defined on \([-3/2,m/2]\) 
and is a solution of (\ref{1bis}).
\par
\bigskip
\noindent
{\bf Remark}: This method explictly requires the ensuing solution to be bounded.
It is possible that the a-priori bound \(M\) is not chosen big enough in the first 
place, i.e., that for some
integer \(m>0\) and for some \(t>0\),
\(
 |x_{m/2}(t)| >M
\), 
in which case, we must reconstruct the sequence
of function (\(x_{m/2}\)) with a 
larger constant \(M\), leading to a smaller value of \(\varepsilon_0\).


\par\vskip10pt\noindent
{\sc Proof of Proposition \ref{p1}}: 
We first need to show that the map $\o_z$ takes its range in 
$\ccol_0([0,1/2])$. Let \(w\in \ccol_0([0,1/2])\). Since
\(
|w(t)|_{0}\leq K_{2},
\)
by choosing 
\[
|\varepsilon|\leq \frac{1}{2K_{2}},
\] 
for all \(0\leq s\leq 1/2\),
we have 
\[
-3/2\leq s-1+ \varepsilon w(s) \leq 0,
\]
and 
therefore the integrant in the right hand side of (\ref{oo}) is well defined.
Also, we clearly have \(\o_{z}(w)(0)=z(0)\)
and
\[ 
|\o_{z}(w)(t)| \leq |z(0)| + \tau \int_{0}^{1/2} 
|F\biggl(
z(s-1+ \varepsilon w(s))
\biggr)
)|ds \leq M+\tau K_{0}/2=K_{2}.
\]
Furthermore, 
we have

\begin{eqnarray*}
  \frac{d}{dt} \biggl(\o_z(w)(t)\biggr)
= \tau 
F\biggl(z(t-1+\varepsilon w(t))\biggr)
\end{eqnarray*}
 and therefore, for all \(0\leq t\leq 1/2\), 

\begin{eqnarray*}
  \biggl| \frac{d}{dt} \biggl(\o_z(w)(t)\biggr)
 \biggr|
\leq \tau 
\biggl|
  F\biggl(z(t-1+\varepsilon w(t))\biggr)
\biggr| \leq \tau K_0\leq \ell
,
\end{eqnarray*}
%
%
%
meaning that  $\o_z$ leaves  $\ccol_0([0,1/2])$ 
invariant and the operator is thus well defined. 
To show that the operator is a contraction,
let $w_1,w_2 \in \ccol_{0}([0,1/2])$. 
We have
\[
\o_x(w_1)(t)-\o_x(w_2)(t)= 
\tau\int_{0}^t\biggl(
G_{w_1}(s)-G_{w_2}(s)
\biggr)ds
\]
where
\[
  G_{w_1}(s)=F\biggl(z(s-1+\varepsilon w_1(s))\biggr)
,\ \ \   G_{w_2}(s)=F\biggl(z(s-1+\varepsilon w_2(s))\biggr)
.\]
Thanks to the Mean Value Theorem, 
we have
\[
  \biggl| 
G_{w_1}(s)-G_{w_2}(s)
  \biggr|_{0}\leq \varepsilon \ell
\sup_{|\xi|\leq M}|F'(\xi)|  \cdot |w_1(t) -w_2(t)|_{0}.
\]
Since 
\[
  |w_1(t) -w_2(t)|_{0} \leq \int_0^{1/2} |w_1^\prime (s)-w_2^\prime(s)|ds \leq \frac{1}{2} 
  {\bf d}(w_1,w_2),
\]
it follows that 
\begin{eqnarray}\label{recfar}
  \biggl| 
G_{w_1}(s)-G_{w_2}(s)
  \biggr|_{0}
  \leq \varepsilon \ell \frac{K_1}{2} {\bf d}(w_1,w_2).
\end{eqnarray}
Therefore 
\[
  {\bf d}\biggl(\o_z(w_1),\o_z(w_2)\biggr) \  
  \leq \ \varepsilon \tau \ell\frac{K_1}{2} {\bf d}(w_1,w_2),
\]
and by choosing 
\[
0<\varepsilon_0\leq \min\{ \frac{1}{2K_{2}}, \frac{1}{\tau \ell K_1} \},
\]
the proof of Proposition 1 is completed.

\subsection{Adding  friction} \label{sec:friction}
We now are in position to understand the dynamics of the full system (\ref{1})
when the friction is taken into consideration.
Recall that the system to be considered satisfies
\begin{equation}\label{withfriction}
\dot x = -a \tau x +\tau F[x(t-1+\varepsilon x(t))], 
\end{equation}
with initial condition $x(t)=x_0(t)$ defined on 
\([-3/2,0]\).
We adapt our previous approach by using 
the method of integrating factor. From (\ref{withfriction}) we write
\begin{eqnarray*}
\frac{d}{dt}\biggl( x(t)e^{a \tau t}\biggr)
=\tau e^{a\tau t} F[x(t-1+\varepsilon x(t))].
\end{eqnarray*}
By integrating both sides of the former equation we get
\begin{eqnarray}\label{wf}
 x(t) = x(0) e^{-a\tau t}+
\tau \int_0^t e^{a\tau (s-t)} F[x(s-1+\varepsilon x(s))]ds.
\end{eqnarray}
To this point our strategy will follow the same as in the former section. 
In what follows, the notation used are the same as before. 
Let \(M>1\),
 $\ell \geq \tau K_0+a\tau M $ and $z \in \ccol([-3/2, 0])$ satisfying (\ref{aboutMM}).
We define the following operator
\[
\jj_z: \ccol[(0,1/2)] \rightarrow  \ccol[(0,1/2)], \ \omega \mapsto \jj_z(\omega)
\]
where
\[
\jj_z(\omega)(t)= z(0) e^{-a\tau t}+
\tau \int_0^t e^{a\tau (s-t)} F[z(s-1+\varepsilon \omega (s))]ds.
\]
The proof of the folllowing proposition is similar to that of Proposition 1 and is left to the reader.
\begin{proposition}\label{p2}
{\sl 
There exists $\varepsilon_0>0$ such that for all
$0\leq\varepsilon\leq \varepsilon_0$,
the operator $\jj_z$ is a contraction. More precisely
for all $\omega_1, \ \omega_2 \in  \ccol_0([0,1/2])$, 
\[
{\bf d}(\jj_z(\omega_1),\jj_z(\omega_2))
\leq \frac{1}{2} {\bf d} (\omega_1,\omega_2).
\]
}
\end{proposition}
As a consequence of Proposition 2, we can construct the Half-Stroboscopic map
\[
\Psi_\varepsilon: \{z\in \ccol([-3/2, 0]) \ | \ \|z\|\leq M \}\rightarrow  \ccol([0,1/2]), \ \ z\mapsto \Psi_\varepsilon(z)
\]
where
$\Psi_\varepsilon(z)$ is the unique fixed point of $\jj_z$, i.e., 
\[
{\Psi}_\varepsilon(z)(t)=z(0)e^{-a\tau t} 
+ \tau \int_{0}^{t} e^{a\tau(s-t)} 
F\biggl(z(s-1+\varepsilon { \Psi}_\varepsilon(z)(s))\biggr)ds, 
\]
and we retrieve the solution of (\ref{withfriction}) using the 
same construction as in Section 3.1.

\section{A numerical method}
%
We now combine the contracting operator 
of the previous section with an interpolating operator
and obtain a finit numerical scheme. 
We will see that the interpolated operator is still a contraction,
arbitrarily close to  the former in the $C^0$ topology.
The  method is presented, for sake of simplicity,
in the friction free case (i.e., 
$a=0$), but is also valid in the case the system admits a friction.

\subsection{The Lagrange-Chebychev interpolation}

Let $q>1$ be an integer and \(L >0\). We denote by   
 ${\bP}_q[t]$ 
 the subset of polynomial functions of degree
less than $q-1$. 
%
We define
the  Lagrange Chebyschev interpolating operator
\[
{\bf {\cal L}}_{q}: \ccoL([0,1/2]) \rightarrow {\bP}_q[t], 
\ h \mapsto {\cal L}_q(h)
\]
where
\[
{\cal L}_q(h)(t)= P_q({\hat h})(4t-1), \ \ \ t\in [0,1/2], \ \ 
P_q({\hat h})(u)=\sum_{j=0}^{q-1} c_j T_j(u),
\]
where 
\[
{\hat h}(u)=h((u+1)/4), \ \ -1\leq u\leq 1,
\]
the $T_j$'s being the Chebyshev polynomial (of the first type) i.e., 
\[
T_j(u) = \cos(j \arccos(u)), \ j=0,\ldots ,q-1, \ \ \ u \in [-1,1],
\]
\[
c_j = \frac{2}{q}\sum_{k=0}^{q-1} {\hat h}(u_k)T_j(u_k), \ \ j >0, \ 
\ c_0 = \frac{1}{q}\sum_{k=0}^{q-1}
{\hat h}(u_k),
\]
where the $u_k$'s are the Chebyshev node on $[-1,1]$, i.e., 
\[
u_k= \cos(\frac{2k+1}{2q}\pi), \ \ \ k=0,\ldots, q-1.
\] 
See \cite{mh6} for more details. 

 \bigskip
 The operator ${\cal L}_q$ is linear and bounded (see below).
Moreover, for all Lipschitz function \(h\), ${\cal L}_q(h)$
 converges uniformly to $h$ on $[0,1/2]$ 
as $q$ tends to $\infty$. More precisely we state the following lemma.
\begin{lemma}
  Let 
  \(
    {\bf g}
    \in C([0,1/2], {\bR})
  \) be differentiable on \((0,1/2)\).
Then
\[
  \sup_{0\leq t \leq 1/2}|
  {\cal L}_q({\bf g})(t)- {\bf g}(t) 
|
  \leq\frac{(1+\mu_q)}{4q} |{\bf g}|_1,
\]
 where
\[
\mu_q= \frac{1}{\pi}\sum_{j=0}^{q-1} \cot\biggl(\frac{(j+1/2)\pi}{2q} \biggr)
=\frac{2}{\pi} \log(q) +0.9625 + {\cal O}(1/q).
\]
\end{lemma}


\bigskip
\noindent
This lemma is a direct consequence of
Jackson's Theorem 
and its 
Corollary 6.14A in \cite{mh6}. Both results are formulated for 
a continuous  function \({\hat {\bf g}}\) defined on \([-1,1]\) with 
modulus of continuity 
\[
{\bf m}(\delta) = \sup_{|x_1-x_2|\leq \delta} 
{|{\hat {\bf g}}(x_1) -{\hat {\bf g}}(x_2)|}.
\]
Corollary 6.14A in \cite{mh6} states that 
\[
  \sup_{|t|\leq 1}\biggl|
{P}_q({\hat {\bf g}})- {\hat {\bf g}} 
\biggr|
\leq {\bf m}(1/q) (1+\mu_q).
\]
After a linear rescaling, 
one  extends these results 
for interpolation on the interval \([0,1/2]\)
and 
in the present case, 
the function \({\hat {\bf g}}\) defined by
\[
  {\hat {\bf g}}(u)=  {\bf g}((u+1)/4) 
\]
satisfies
 \begin{eqnarray}\label{bdbis}
   \sup_{-1<t<1} |{\bf {\hat g}}^\prime(u)|
   = |{\bf g}|_1/4
  \end{eqnarray}
and therefore
\[
  {\bf m}(1/q)\leq  \frac{|{\bf g}|_1}{4q}, 
\]
and the lemma follows. 
Finally, observe that the estimate 
\begin{eqnarray}\label{mu}
\mu_q =\frac{2}{\pi} \log(q) +0.9625 + {\cal O}(1/q)
\end{eqnarray}
is also given in \cite{mh6}, 
which implies that 
\[
|
{\cal L}_q({\bf g})- {\bf g}
|_{0}\to 0 \ \  {\mathrm {as}}\  q  \to \infty.
\]
%
 Observe that the upper bound for 
  \( 
  |{\cal L}_q({\bf g}) -{\bf g}|_0\)
given in Lemma 1 is often 
far from being 
optimal.  For instance,
following 
\cite{mh6},
  if \({\bf g}\) is analytic, there exists 
\({\tilde M}>0\) and \({\tilde r}>1\) such that 
\[
  |{\cal L}_q({\bf g})-{\bf g}|_{0} \leq {\tilde M}/{\tilde r}^q.
\]
In the latter estimate,
both constants  \({\tilde M}\) and \({\tilde r}\) 
depend upon the function \(z\). However, in the case 
  \({\bf g} \in {\bf C}^{0,L}([0,1/2])\),
thanks fto Lemma 1, we have the following estimate
  \[
    |{\cal L}_q({\bf g})-{\bf g}|_{0} \leq L\frac{1+\mu_q}{4q}  
,
\]
and the upper bound only depends upon the constant \(L\).
Indeed, we state the following corollary.

\begin{cor}\label{1111}
    Let \(L>0\).
Then there exists $q_0=q_0(L)>1$ such that 
for all $q\geq q_0$ and for all 
  \( 
  {\bf g} \in {\bf C}^{0,L}([0,1/2])\),
\[
  |
  {\cal L}_q({\bf g})-{\bf g}
|_{0}  
  \leq 
  K_0 
.
\]
\end{cor}

\subsection{The reduced  operator}
The limitation we have with the contraction operator introduced
in the former section is that we have to work in
an infinite dimension space.
To overcome this difficulty
we replace the contraction operator by 
the so called {\it reduced operator}. The latter
is an approximation of the former. We show that 
the reduced operator is also a contraction.

\bigskip

Let \(M>1\), 
\(\ell \geq 2\tau K_0\) and  
\(z \in \ccol([-3/2,0])\) satisfying (\ref{1bis}).
We also  assume \(z\) to admits a bounded 
  second derivative on \((-3/2,0)\).
  As a consequence, there exists positive constants 
\(l_1\) and \(l_2\)
such that the function 
\[
  {\bf H}: [-3/2,0] \rightarrow {\bR}, \ u\mapsto  F(z(u))
\]
satisfies
\begin{eqnarray}
  \label{2derf}
  \biggl|
  {\bf H}^{\prime} (u)
  \biggr| \leq  l_1, \ \ \  
  \biggl|
  {\bf H}^{\prime\prime}(u) 
  \biggr| \leq  l_2, \ \ \ {\mathrm {for}}\ {\mathrm{all}} \ -3/2<u<0. 
\end{eqnarray}
We  define 
 \[
\o_{z,q}:
{\bP}_{q}[t] \cap \ccol_0([0,1/2])\rightarrow {\bP}_{q}[t] \cap
\ccol_0([0,1/2])
,
\ {\bf f} \mapsto \o_{z,q} ({\bf f})
 \]
%
where \(\ccol_0([0,1/2])\) is defined in (\ref{k22}) 
with (unlike in the former section)
 \[
 K_{2}=M+\tau K_{0},
 \]
 where
\begin{eqnarray}\label{last}
\o_{z,q}({\bf f})(t) & = & {\bf f}(0) +
\tau\int_{0}^t {\cal L}_{q-1}
  \biggl(G_{\bf f}(s)
  \biggr)ds.
\end{eqnarray}
and where
\[
  G_{\bf f}(t)= F(z(t-1+\varepsilon {\bf f}(t))). 
\]
We state the following proposition.

\begin{proposition}\label{p3}
{\sl 
  There exists \(q_1\geq q_0\),
such that for all
  \(
  0\leq|\varepsilon|\leq \varepsilon_0
  \), 
  and for all \(q\geq q_1\),
  $\o_{z,q}$ is a contraction. More precisely
for all 
  \(
  \omega_1, \ \omega_2 \in  {\bP}_{q,0}[t] \cap   \ccol_{0}([0,1/2])
  \), 
\[
  {\bf d}
  \biggl(
  \o_{z,q}(\omega_1),\o_{z,q}(\omega_2)
  \biggr)
 \leq \frac{1}{2} 
  {\bf d}(w_1,w_2)
.
\]
}
\end{proposition}

\noindent
Before proving this proposition, we need the following lemma.
\begin{lemma}
Let 
  \[
    {\bf H} : [-3/2,0] \rightarrow {\bR}, \ t \mapsto {\bf H}(t)
  \]
  be a continuous, twice differentiable map on \( (-3/2,0)\). 
  Assume that 
  \[
   \biggl|
  {\bf H}^{\prime} (u)
  \biggr| \leq  k_1, \ \ \  
  \biggl|
  {\bf H}^{\prime\prime}(u) 
  \biggr| \leq  k_2, \ \ \ {\mathrm {for}}\ {\mathrm{all}} \ -3/2<u<0,
  \]
  for some \(k_1>0\) and \(k_2>0\).
  Let \(L>0\), \(0\leq |\varepsilon|<\varepsilon_0\),  \(w_1, w_2 \in {\bf C}^{0,L}_0([0,1/2])\) satisfying 
  \[
    \max\{ |w_1(t)|_{0}, |w_2(t)|_{0}\} \leq M 
  .\]
  Then, 
  \[
  \biggl|
  \biggl({\bf H}\circ W_1 -{\bf H}\circ W_2
  \biggr)^{\prime}
  \biggr|_{1} 
  \leq
  \varepsilon \biggl(\frac{(1+\varepsilon L)k_2}{2} +k_1\biggr) 
  {\bf d}(w_1,w_2),
  \]
where
  \[
    W_1(t)=t-1+\varepsilon w_1(t) , \ \ \  \\ \ W_2(t)=t-1+\varepsilon w_2(t)
 . \]
\end{lemma}
\noindent
{\bf Proof}: 
By choosing 
\(|\varepsilon| \leq \varepsilon_0\), 
for all \(0\leq s\leq 1/2\), we have
\[
  -3/2\leq W_i(s)=s-1+\varepsilon w_i(s)\leq 0,\ \ \ \ \ \ \ \ \ \ \ \ \ \ i=1,2
\]
and therefore both \({\bf H}\circ W_1\) and \({\bf H}\circ W_2\) are 
well defined.
We write
\begin{eqnarray*}
  \biggl(
    {\bf H}\circ W_1 -
    {\bf H} \circ W_2 
  \biggr)^{\prime}(t)
   & = & 
  W^\prime_1(t) 
   {\bf H}^\prime(W_1(t)) - 
  W^\prime_2(t) 
   {\bf H}^\prime(W_2(t))\\ 
  & = & 
  W^\prime_1(t)
  \biggl( 
   {\bf H}^\prime(W_1(t)) - 
   {\bf H}^\prime(W_2(t))\biggr)
   \\
  & + & 
  \biggl(
  W^\prime_1(t)-W^\prime_2(t)\biggr) 
   {\bf H}^\prime(W_2(t))
   .
\end{eqnarray*}
By the Mean Value Theorem,  for all \(0<t<1/2\), we have
\begin{eqnarray*}
  \biggl|
 \biggl(
    {\bf H}\circ W_1 -
    {\bf H} \circ W_2 
  \biggr)^{\prime}(t)
\biggr|
  & \leq & 
   \sup_{-3/2< u< 0}|{\bf H}^{\prime\prime}(u)| 
  \biggl|W_1(t) -W_2(t)\biggr|
 \biggr|W^\prime_1(t)\biggr| 
\\
  & + & 
   \sup_{-3/2< u< 0}|{\bf H}^{\prime}(u)| 
  \biggl|
  W^\prime_1(t) 
  -
  W^\prime_2(t)\biggr| 
\end{eqnarray*}
and 
since
\[
  |W_1^\prime(t)| \leq 1+\varepsilon L, \ 
\]
we have
\begin{eqnarray}\nonumber 
  \biggl|
 \biggl(
    {\bf H}\circ W_1 -
    {\bf H} \circ W_2 
  \biggr)^{\prime}(t)
\biggr|
&  \leq & 
   k_2
  (1+\varepsilon L)
  \biggl|W_1(t) -W_2(t)\biggr|
\\ \label{wh}
& + & 
 k_1 
  \biggl|
  W^\prime_1(t) 
  -
  W^\prime_2(t)\biggr|. 
\end{eqnarray}
Furthermore, for all \(0< t<1/2\)
\begin{eqnarray}\label{m11}
  \biggl|
  W^\prime_1(t) 
  -
  W^\prime_2(t)\biggr|=
\varepsilon
  \biggl| 
  w^\prime_1(t)-w_2^\prime(t)
  \biggr|\leq \varepsilon{\bf d}(w_1,w_2) 
,
\end{eqnarray}
and
\begin{eqnarray}\nonumber
  \biggl|
  W_1(t) 
  -
  W_2(t)\biggr|=
\varepsilon
  \biggl| 
  w_1(t)-w_2(t)
  \biggr| & \leq & \varepsilon \int_0^{1/2} 
  |
  w^\prime_1(s)-w_2^\prime(s)
  |
  ds 
  \\ \label{m22}
  & \leq & 
  \frac{\varepsilon}{2} {\bf d}(w_{1},w_{2})
.
\end{eqnarray}
Thanks to (\ref{wh}), (\ref{m11}) and ({\ref{m22}) it follows
that for all \(0< t< 1/2\), 
\begin{eqnarray*}
  \biggl|
 \biggl(
    {\bf H} \circ W_1 -
    {\bf H}\circ W_2 
  \biggr)^{\prime}(t)
\biggr|
  \leq 
   \varepsilon \frac{k_2}{2}
  (1+\varepsilon L)
  {\bf d}(w_1,w_2)
   +  
 \varepsilon k_1 
  {\bf d}(w_1,w_2)
,
\end{eqnarray*}
ending the proof of the lemma.

%
%

\bigskip
\noindent
{\sc Proof of Proposition 3}:
We first verify that $\o_{z,q}$ leaves 
${\bP}_{q}[t] \cap \ccol_0([0,1/2])$ invariant. 
For all $0\leq t\leq 1/2$  and for all $\omega \in \ccol_0([0,1/2]$
from its definition 
$\o_{z,q}(w)$ is a polynomial of degree equal or less than $q-1$ and 
$\o_{z,q}(w)(0)=z(0)$.
Again,
by choosing 
\(|\varepsilon| \leq \varepsilon_0 \),  
for all \(0\leq s\leq 1/2\) we have
\[
-3/2\leq s-1+\varepsilon \omega(s)\leq 0,
\]
%
and thus  the integrant in the right hand side 
of (\ref{last}) is well defined.
Observe that function
\[
s\mapsto G_\omega(s)= F\biggl(
    z(s-1+ \varepsilon \omega(s))
    \biggr) 
\]
is 
differentiable
and for all \(0<t<1/2\),
\[
\biggl|
G^\prime_\omega(s)
\biggr| \leq L=K_1(\ell(1+\varepsilon \ell),
\]
and bounded by \(K_0\). 
By taking \(q\geq q_{0}+1\) given in Corollary 1, 
for all \(0\leq t\leq 1/2\),
we have
\begin{eqnarray*}
\begin{array}{rcl}
\biggl|
  \o_{z,q}(\omega)(t)
  \biggr| & \leq  & |z(0)| + \tau \displaystyle\int_{0}^{1/2} 
\biggl|{\cal L}_{q-1} \biggl( G_\omega(s) \biggr)-G_\omega(s)
\biggr|ds
\\ &  &  \\
& +  & 
 \tau \displaystyle\int_{0}^{1/2} 
 \biggl|
  G_\omega(s)
  \biggr| 
  ds 
\\ \ & & \\
& \leq & M + \displaystyle\frac{\tau}{2}\cdot(2K_{0}) = K_{2}.
\end{array}
\end{eqnarray*}
Furthermore, thanks again to Corollary 1, for all \(0<t<1/2\), we have
\begin{eqnarray*}
  \biggl| 
  \biggl(\o_{z,q}(w)\biggr)^{\prime}(t)
  \biggr|  
  & \leq   &  
\tau \Biggl| 
  {\cal L}_{q-1}
\biggl(G_{w}(t)
\biggr)
\biggr)
\Biggr| \\ 
& \leq &   
  \tau 
\Biggl| 
  {\cal L}_{q-1}
\biggl(
G_{w}(t)\biggr) - G_{w}(t) 
\Biggr|
+  
\tau 
\biggl| 
G_{w}(t) 
\biggr|
\\ 
& \leq & 2 \tau K_0 \leq \ell, 
\end{eqnarray*}
meaning that  $\o_z$ leaves  
$\ccol_0([0,1/2])$ invariant. 
We now show that the operator is a contraction.
Let $w_1,w_2 \in \ccol_{0}([0,1/2])$. A straightforward computation 
gives
\begin{eqnarray}\label{fj}
  \biggl(
  \o_{z,q}(w_1)(t)-\o_{z,q}(w_2)
\biggr)^{\prime}(t)
  = 
\tau
  {\cal L}_{q-1}\biggl(G_{w_1}
  -G_{w_2}\biggr)(t).
\end{eqnarray}
Thanks to Lemma 2 and since
\(G_{w_1}(0)=G_{w_2}(0)\),
for all \(0< t<  1/2\),
\begin{eqnarray}\label{almost}
\Biggl|
\biggl(
G_{w_1}-G_{w_2}
\biggr)^{\prime}(t)
\Biggr| \leq 
\biggl|
G_{w_1}-G_{w_2}
\biggr|_1 = 
{\bf d}(G_{w_1}, G_{w_2}) 
\leq
\varepsilon {\bf L} {\bf d}(w_1,w_2)
,
\end{eqnarray}
where
\[
  {\bf L}=\frac{(1+\varepsilon\ell)l_2}{2}+l_1.
\]
Thanks to Lemma 1, (\ref{recfar}) and (\ref{almost}), 
for all \(0\leq t \leq 1/2\), 
\begin{eqnarray*}
\Biggl|
{\cal L}_{q-1}\biggl(
G_{w_1}-G_{w_2}
 \biggr)(t)
\Biggr|
 & \leq &  
\biggl|
  {\cal L}_{q-1}\biggl(
G_{w_1}-G_{w_2}
 \biggr)(t)-
  \biggl(
  G_{w_1}(t)-G_{w_2}(t)
  \biggr)\biggr|
  \\ &  + & 
\biggl|
  G_{w_1}(t)-G_{w_2}(t) 
  \biggr|
\\
  & \leq &  \biggl(\frac{1+\mu_q}{4q}
\biggr){\bf d}(G_{w_1}, G_{w_2})
  + \varepsilon \ell\frac{K_1}{2} {\bf d}(w_1,w_2)
  \\ & \leq &
\varepsilon\biggl(\frac{1+\mu_q}{4q}
   {\bf L} +\ell\frac{K_1}{2}
\biggr)
  {\bf d}(w_1,w_2)
.
\end{eqnarray*}
\noindent
We now take \(q_{1}\) sufficiently large such that. for all \(q\geq q_{1}\) satifying
\[
\frac{1+\mu_{q}}{4q} L \leq \ell K_{1}/2,
\]
and from above we get
\begin{eqnarray*}
\Biggl|
{\cal L}_{q-1}\biggl(
G_{w_1}-G_{w_2}
 \biggr)
\Biggr|_{0} \leq 
 \varepsilon \ell K_{1} {\bf d}(w_1,w_2)
.
\end{eqnarray*}
Finally with (\ref{fj}), it follows that
\[
{\bf d}\biggl( \o_{z,q}(w_1),\o_{z,q}(w_2)\biggr) \leq 
\tau \int_{0}^{1/2}
\Biggl|
{\cal L}_{q-1}\biggl(
G_{w_1}-G_{w_2}
 \biggr)
\Biggr|(s) ds
\leq 
\varepsilon \tau \ell\frac{K_{1}}{2} {\bf d}(w_1,w_2)
.\]
By choosing \(0\leq |\varepsilon|\leq \varepsilon_{0}\),  Proposition 3 is proved. 
As a consequence, 
for 
each $z \in \ccol([-3/2 ,0])$, the operator $\o_{z,q}$ admits a unique fixed point 
and one deduces the same construction as in section 3.2 replacing 
$\o_{x}$ by $\o_{x,q}$.

\subsection{Constructing the orbit}
Fix $\nu$ a small positive number representing the 
tolerance of our computation. 
Take $\ell>0$ as above and $x_0 \in \ccol[-3/2, 0]$.
Our goal is now 
 to compute 
$\Psi_\varepsilon(x_0)$ with an arbitrary accuracy, more precisely
we aim to 
compute 
$y_0=\Psi_\varepsilon(x_0)$ and more precisely to 
find a function
${\tilde y}_{0}$
such that
\begin{eqnarray}\label{final}
    {\bf d}\biggl(\Psi_\varepsilon(x_0),{\tilde y}_0\biggr)
\leq \nu.
\end{eqnarray}
Let $q >1$ 
be an integer such that
for all $\omega \in \ccol([0,1/2])$ and independently from the choice of \(\varepsilon\), 
\begin{eqnarray}\label{l-e1}
 \biggl| {\cal L}_{q-1} \biggl(G_\omega\biggr) -G_\omega
    \biggr |_0  
\leq \nu/(4\tau).
\end{eqnarray}
The existence of such integer $q$ is guaranteed by Lemma 1.
From Eq. (\ref{l-e1}), for all $\omega\in \ccol([0,1/2])$ we have
\begin{eqnarray}\label{l-e2}
{\bf d}\biggl(\o_{x_0}( \omega ) ,\o_{x_0,q}(\omega)\biggr) \leq \nu/4.
\end{eqnarray}

Let ${\bf f}_0 \in \ccol([0,1/2])$ and construct the
sequence of functions 
\[
{\bf f}_n \in \ccol([0,1/2]),\ \ {\mathrm{such}} \ {\mathrm{that}} \ \ 
{\bf f}_{n+1}(t) =\o_{x_0,q} ({\bf f}_n)(t).
\]
Thanks to Proposition 3, there exists an integer ${\bf m} \geq 1$ 
such that 
\begin{eqnarray}\label{f-e3}
{\bf d}\bigg(
\o_{x_0,q}({\bf f}_n) - {\bf f}_n)
\biggr) \leq \nu/4, \ \ \ 
\forall n \geq {\bf m}.
\end{eqnarray}
We now write ${\tilde y}_0={\bf f}_{\bf m}$.
By definition we have
\[
\o_{x_0}(\Psi_\varepsilon(x_0))= \Psi_\varepsilon(x_0)
\]
and we also have 
\begin{eqnarray}\label{f2}
 {\bf d}\biggl(\Psi_\varepsilon(x_0),{\bf f}_{\bf m}\bigg)  
 & = &  {\bf d}\biggl(  \o_{x_0}(\Psi_\varepsilon(x_0)),{\bf f}_{\bf m} \biggr) 
\\ \nonumber 
 & \ & 
\\ \nonumber 
& \leq &  
{\bf d}\biggl( \o_{x_0}(\Psi_\varepsilon(x_0)),\o_{x_0}({\bf f}_{\bf m})\biggr)\\ \nonumber
& \ & \\ \nonumber
& + &  {\bf d}\biggl(  \o_{x_0}({\bf f}_{\bf m}),  \o_{x_0,q}({\bf f}_{\bf m}) \biggr) 
 + {\bf d}\biggl(   \o_{x_0,q}({\bf f}_{\bf m}) , {\bf f}_{\bf m} \biggr).
\end{eqnarray}
From Proposition 1 we have
\begin{eqnarray}\label{neee}
{\bf d}\biggl(\o_{x_0}(\Psi_\varepsilon(x_0)), \o_{x_0}({\bf f}_{\bf m}) \biggr)
\leq  \frac{1}{2}{\bf d}\biggl( \Psi_\varepsilon(x_0),  {\bf f}_{\bf m} \biggr).
\end{eqnarray}
Furthermore from (\ref{l-e2}) we have 
\begin{eqnarray}\label{recalll}
{\bf d}\biggl(
 \o_{x_0}({\bf f}_{\bf m}),  \o_{x_0,q}({\bf f}_{\bf m}) 
\biggr) 
\leq \nu/4.
\end{eqnarray}
Finally, thanks to (\ref{f2}), (\ref{f-e3}), (\ref{neee}) and   (\ref{recalll}), 
we have
\begin{eqnarray*}\label{f12}
    {\bf d}\biggl(\Psi_\varepsilon(x_0) , {\bf f}_{\bf m}\biggr)    & \leq  &  
    \frac{1}{2}{\bf d}\biggl( \Psi_\varepsilon(x_0) , {\bf f}_{\bf m} \biggr)
     +\nu/4 +\nu/4
%
   \end{eqnarray*}
and therefore 
\[
{\bf d}\biggl( \Psi_\varepsilon(x_0) ,{\tilde y}_0\biggr)\leq \nu.
\]
This above procedure allows us to construct $\Psi_\varepsilon(x_0)\sim {\tilde y}_0$, we then deduce
\[
x_{1/2} : [-1 +\varepsilon x_0(0) -1/2, 0] \rightarrow {\bR} ,  \ t 
\mapsto x(t+1/2)
\]
where
\[
x(t) = x_0(t) \ \ \ {\mathrm {if}} \ t<0, \ \ x(t) = {\tilde y}_0(t) 
\ \ \ {\mathrm {if}} \ 0\leq t\leq 1/2.
\]
Following the above construction we deduce
\({\tilde y}_{1/2}\) which approximate \(\Psi_\varepsilon(x_{1/2})\) up to a \(\nu\)-tolerance 
and following the same notation as in section 3.2, we can retrieve the 
solution \(x\) on the entire real line.

\section{Numerical Results}
\label{sec:numsim}

\begin{table}
\centering
\small
\begin{tabular}{ccc}
\toprule
System & System Delay ($\tau$) & State-Dependence ($\varepsilon$) \\
\midrule
\multirow{4}{*}{Cubic Ikeda} & \multirow{4}{*}{1.62} & 0  \\
 & & 0.05 \\
 & & 0.15 \\
 & & 0.20  \\
\midrule
\multirow{6}{*}{Mackey-Glass} & \multirow{3}{*}{2} & 0 \\
 & & 0.05  \\
 & & 0.15 \\
\cmidrule[1pt](lr){2-3}
 & \multirow{3}{*}{4} & 0\\
 & & 0.05 \\
 & & 0.10 \\
\bottomrule
\end{tabular}
\caption{Systems and system parameters explored numerically in this study.}
\label{tab:numexp}
\end{table}

Although the results in the preceding sections establish asymptotic convergence of the numerical approximation to the true solution as the degree of the Lagrange Chebyschev interpolating operator is taken to infinity, the behavior of the reduced operator, and the accuracy of the approximation in practice for small degrees is not understood. To empirically investigate this convergence we compared numerically-derived attractors of the cubic Ikeda and Mackey-Glass systems for finite but large $q$, to the derived attractors for small $q$, using several qualitative and quantitative measures of similarity over a range of system and state-dependent delay parameters. In particular, we estimated 1) peak-to-peak maps (PP) over large time intervals to provide a qualitative comparison of system dynamics as the degree of interpolation increases, 2) the correlation dimensions (CD) of delay embeddings of solution trajectories and 3) quantified measures of the geometric and topological similarities of the attractors using persistent homology (PH). 

For each system and parameter choice in Table \ref{tab:numexp} we instantiated ``ground truth'' simulations with random interpolating polynomial coefficients, with $q = 17$ corresponding to degree 17 Lagrange interpolating polynomials defined on [0,1/2]. Thus, the highest accuracy solutions were constructed from concatenations of the derived interpolating polynomials with a total of 34 interpolating nodes per unit interval. For each simulation we then truncated the initial polynomial coefficients to the leading $q = 2,3,4,5,6,7,8,9$ and $10$ coefficients, and simulated coarser approximations of the ground truth solution. All simulations were carried out for 2.1e4 steps, and only the final 2e4 steps were retained to minimize the effect of transients. This yielded solution curves each covering an interval of length 1e4. All simulations used 30 Picard iterations per step.

We first visualize the ground truth systems and illustrated their range of dynamic behaviors using Lissajou plots and peak-to-peak maps.  The Lissajou plots are plots of the 
parametric curve
\[
\{ (x(t), x(t-1)), \ | \ 0 \leq t \leq t_{\text{max}}.\} 
\]
Since time was rescaled by a factor $\tau$, this indeed
represents
\[
\{ y(t), y(t-\tau )), \ | \ 0 \leq t \leq t_{\text{max}}\tau \}. 
\]
Peak-to-peak maps are scatter plots of consecutive local maxima of 1-dimensional solution trajectories. Many choatic systems (e.g., Lorez \cite{DeterministicNonperiodicFlow}) are known to exhibit so-called peak-to-peak dynamics in that the value and time of future local maxima can be accurately predicted from the value of previous local maxima, or peaks. A peak-to-peak map, which relates the values of a 1-dimensional trajectory at consecutive local maxima, represents a reduced order model of the original system which can been used, for instance, in optimal control \cite{PICCARDI2000298}. Here, we use peak-to-peak maps to provide a compact visualization of the system, allowing a by-eye comparison of the approximating-system dynamics to the true system, and a qualitative indication of system stabilization as the number of interpolating nodes is increased. To construct peak-to-peak maps, local maxima were extracted from interpolated solution trajectories over the time interval [0,2.5e3], sampled at a resolution of 100 points per unit interval. 

Restricting first to the high-accuracy ($q=17$) simulations, we observe qualitative changes in dynamic behavior as the state-dependence parameter $\varepsilon$ increases. For example, correlation dimension estimates, peak-to-peak and Lissajou plots suggest chaotic behavior for the Ikeda cubic system with $\varepsilon=0, 0.05$ and $0.2$, while for $\varepsilon = 0.1$ the system appears to be periodic (See Table \ref{tab:quant_metrics_cubic} and Figure \ref{fig:cubic_ground_truth}). We observe a similar bifurcation to simpler dynamics between $\varepsilon = 0.05$ and $0.1$ for the Mackey-Glass with system delay $\tau = 4$ (See Table \ref{tab:quant_metrics_mg_4} and Figure \ref{fig:mg_ground_4}). For Mackey-Glass with system delay $\tau = 2$ there is little change in the apparent and estimated dimension over the range of simulated state-dependent delay parameters,$\varepsilon=0, 0.05$, and $0.15$ (See Table \ref{tab:quant_metrics_mg_2} and \ref{fig:mg_ground_2}).

\begin{figure}[ht]
\begin{center}
\includegraphics[width=1.0\textwidth]{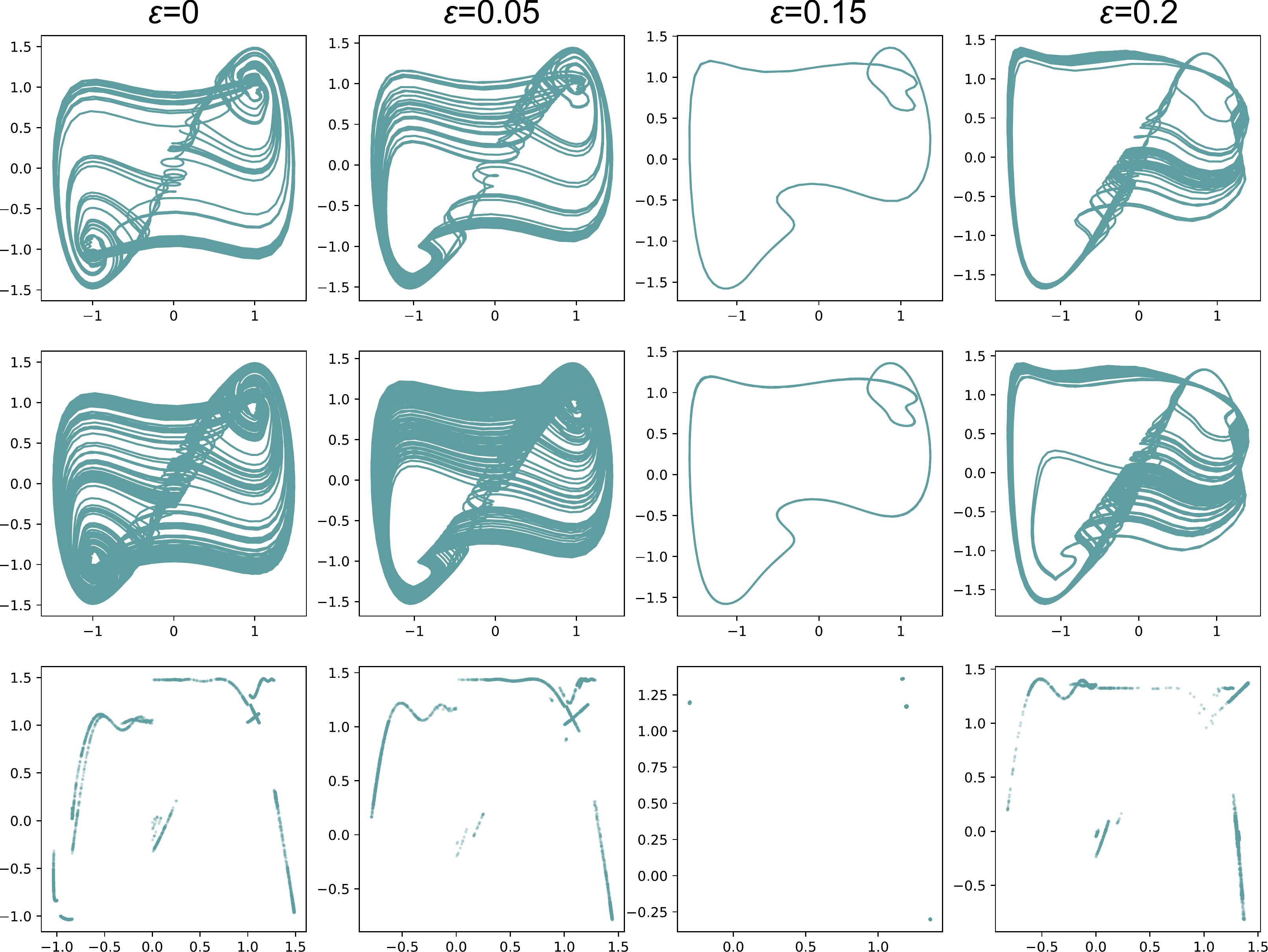}
\caption{Lissajou plots and peak-to-peak maps for the ground-truth ($q=17$) Ikeda cubic system with system delay parameters $\tau=1.62$ and (left to right) $\varepsilon=0, 0.05, 0.15, 0.2$. The top row contains Lissajou plots over the time interval [0,250], the middle row contains Lissajou plots over the time interval [0,1e3], and the bottom row contains peak-to-peak maps for peaks extracted from solutions over the time interval [0,2.5e3].}
\label{fig:cubic_ground_truth}
\end{center}
\end{figure}

\begin{figure}[ht]
\begin{center}
\includegraphics[width=1.0\textwidth]{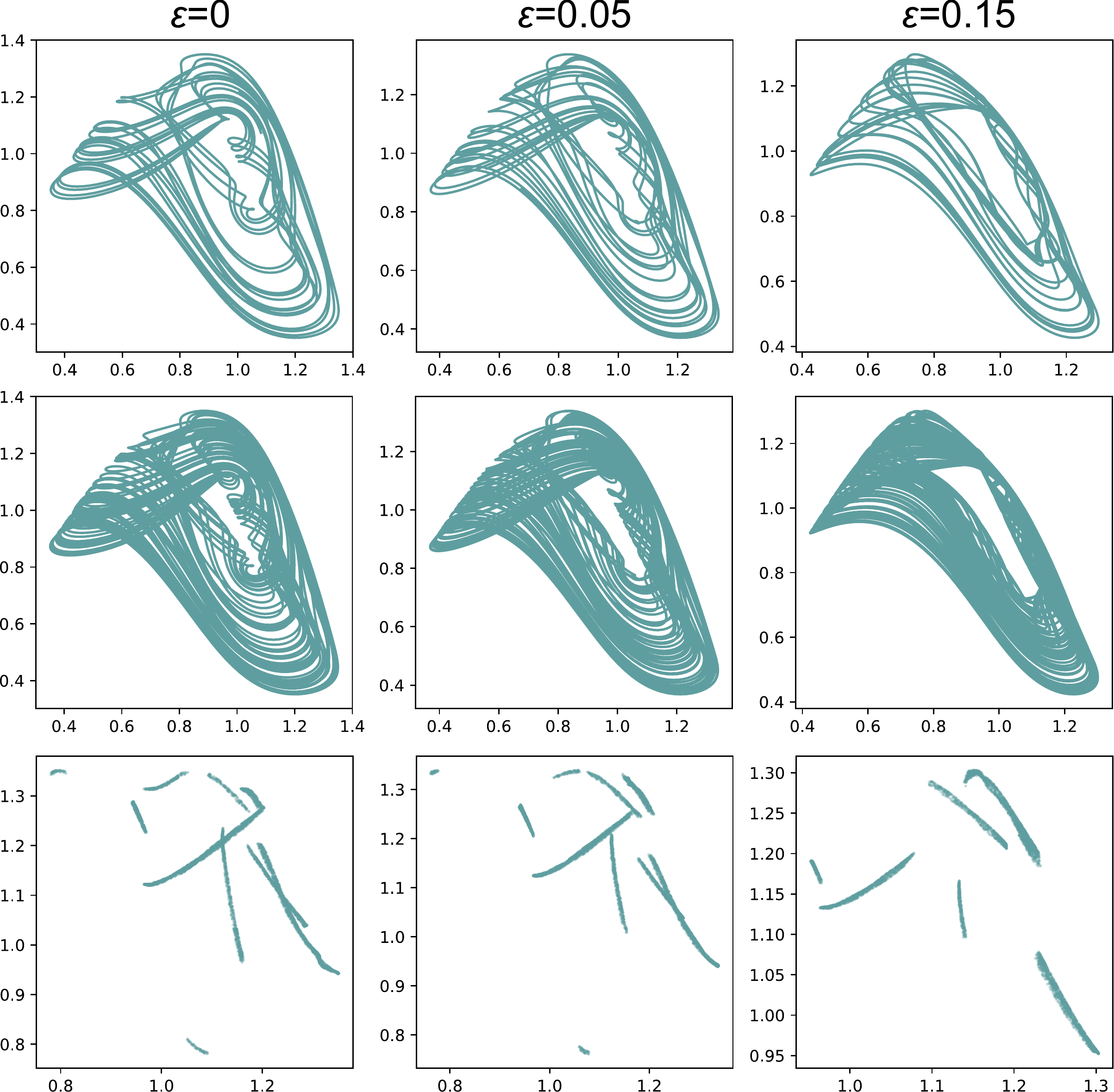}
\caption{Lissajou plots and peak-to-peak maps for the ground-truth ($q=17$) Mackey-Glass system with system delay parameters $\tau=2$ and (left to right) $\varepsilon=0, 0.05, 0.15$. The top row contains Lissajou plots over the time interval [0,50], the middle row contains Lissajou plots over the time interval [0,150], and the bottom row contains peak-to-peak maps for peaks extracted from solutions over the time interval [0,2.5e3].}
\label{fig:mg_ground_2}
\end{center}
\end{figure}

\begin{figure}[ht]
\begin{center}
\includegraphics[width=1.0\textwidth]{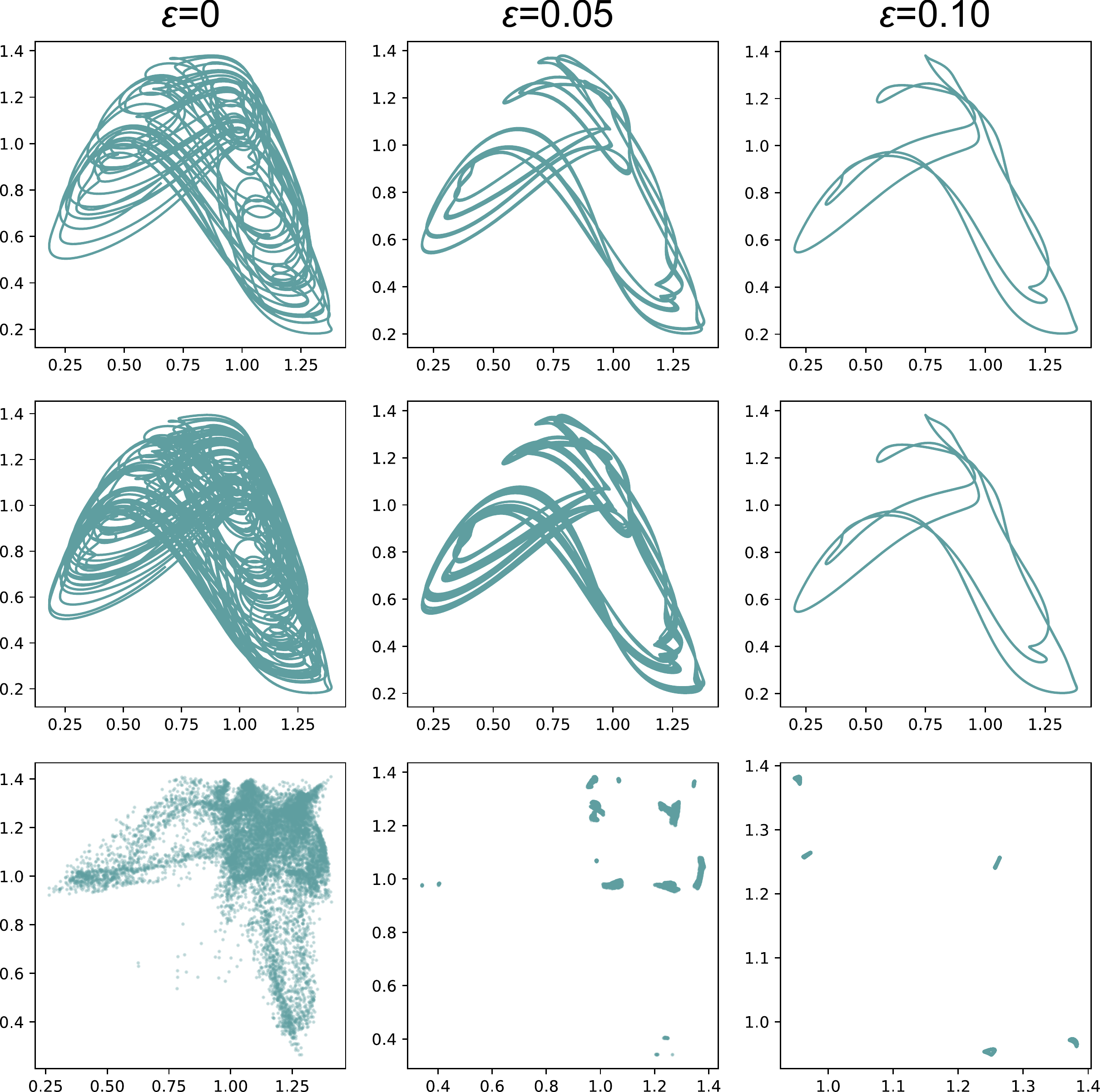}
\caption{Lissajou plots and peak-to-peak maps for the ground-truth ($q=17$) Mackey-Glass system with system delay parameters $\tau=4$ and (left to right) $\varepsilon=0, 0.05, 0.1$. The top row contains Lissajou plots over the time interval [0,50], the middle row contains Lissajou plots over the time interval [0,100], and the bottom row contains peak-to-peak maps for peaks extracted from solutions over the time interval [0,2.5e3].}
\label{fig:mg_ground_4}
\end{center}
\end{figure}

A striking feature that is consistent across these systems and system parameters is the apparent, often dramatic, convergence/stabilization of peak-to-peak dynamics as one increases the number of interpolating nodes, $q$, and a variety of changes in apparent system complexity for small $q$. For example, for the Ikeda system with $\varepsilon=0$ and $q<4$, the peak-to-peak dynamics are markedly different than for the ground truth system, and suggest reduced dynamic complexity. However, for $q \geq 4$ there is little perceptible difference between the approximations and the $q=17$ map (Figure \ref{fig:cubic_162_0_p2p_comps}). 

\begin{figure}[ht]
\begin{center}
\includegraphics[width=1.0\textwidth]{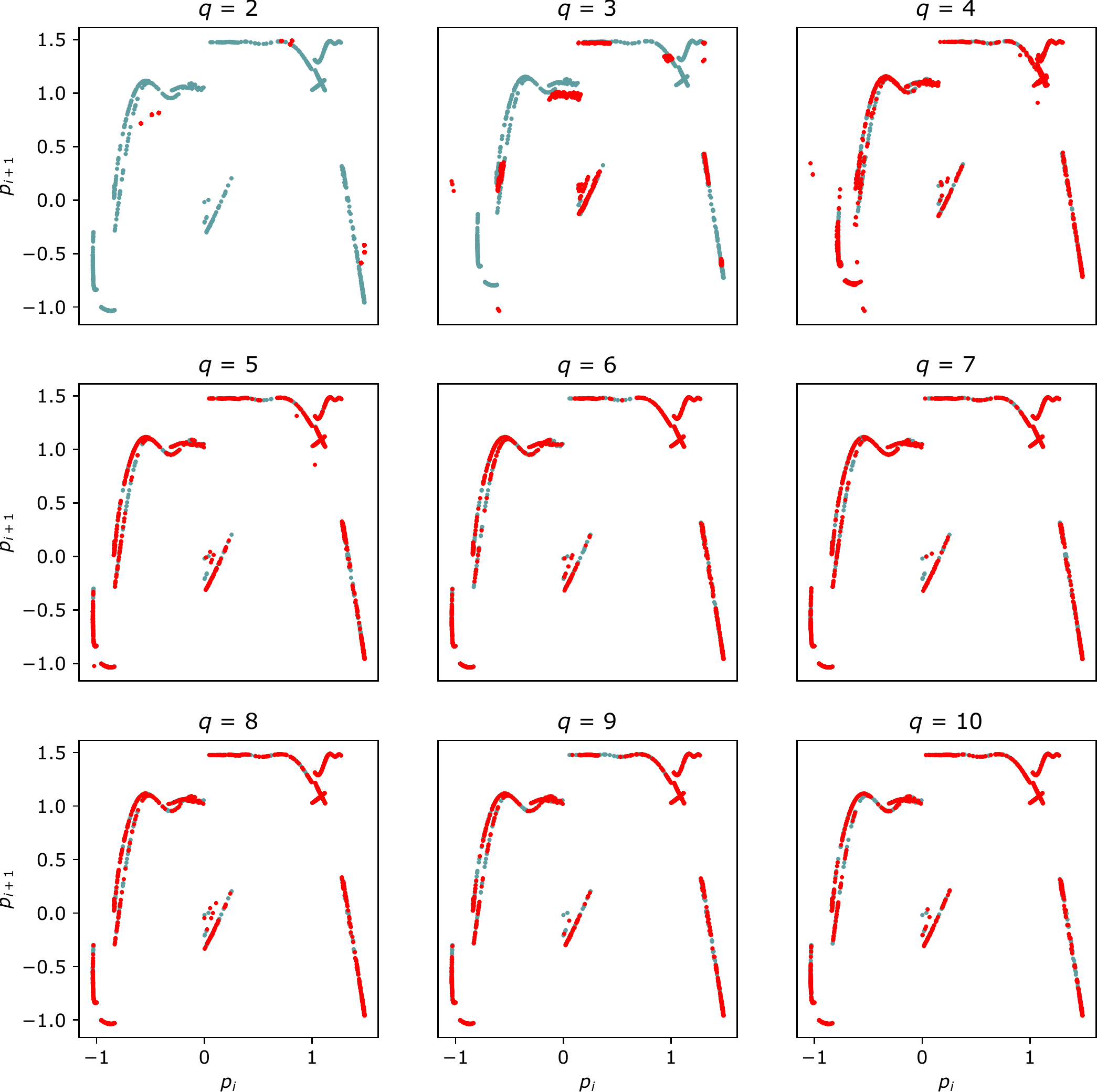}
\caption{Peak-to-peak maps for the ground-truth ($q=17$) Ikeda system with system delay parameters $\tau=1.62$  and  $\varepsilon=0$ (blue dots) compared to approximating system ($q=2,\ldots,10$, red dots). Peaks extracted from solutions over the time interval [0,2.5e3].}
\label{fig:cubic_162_0_p2p_comps}
\end{center}
\end{figure}

In contrast, the true Mackey-Glass system with $\tau=4$ and $\varepsilon = 0.1$ appears to be periodic, while the approximate systems $q=3,4,5,6$ appear to be highly chaotic (Figure \ref{fig:mg_401_p2p_comps}). Interestingly, the dynamics appears to have become periodic by $q=7, 8, 9$ but with extra peak pairs compared to the true system. These spurious peaks appear to be entirely eliminated by $q=10$, at which point there appears to be good agreement with the true system, suggesting the approximation has stabilized.

\begin{figure}[ht]
\begin{center}
\includegraphics[width=1.0\textwidth]{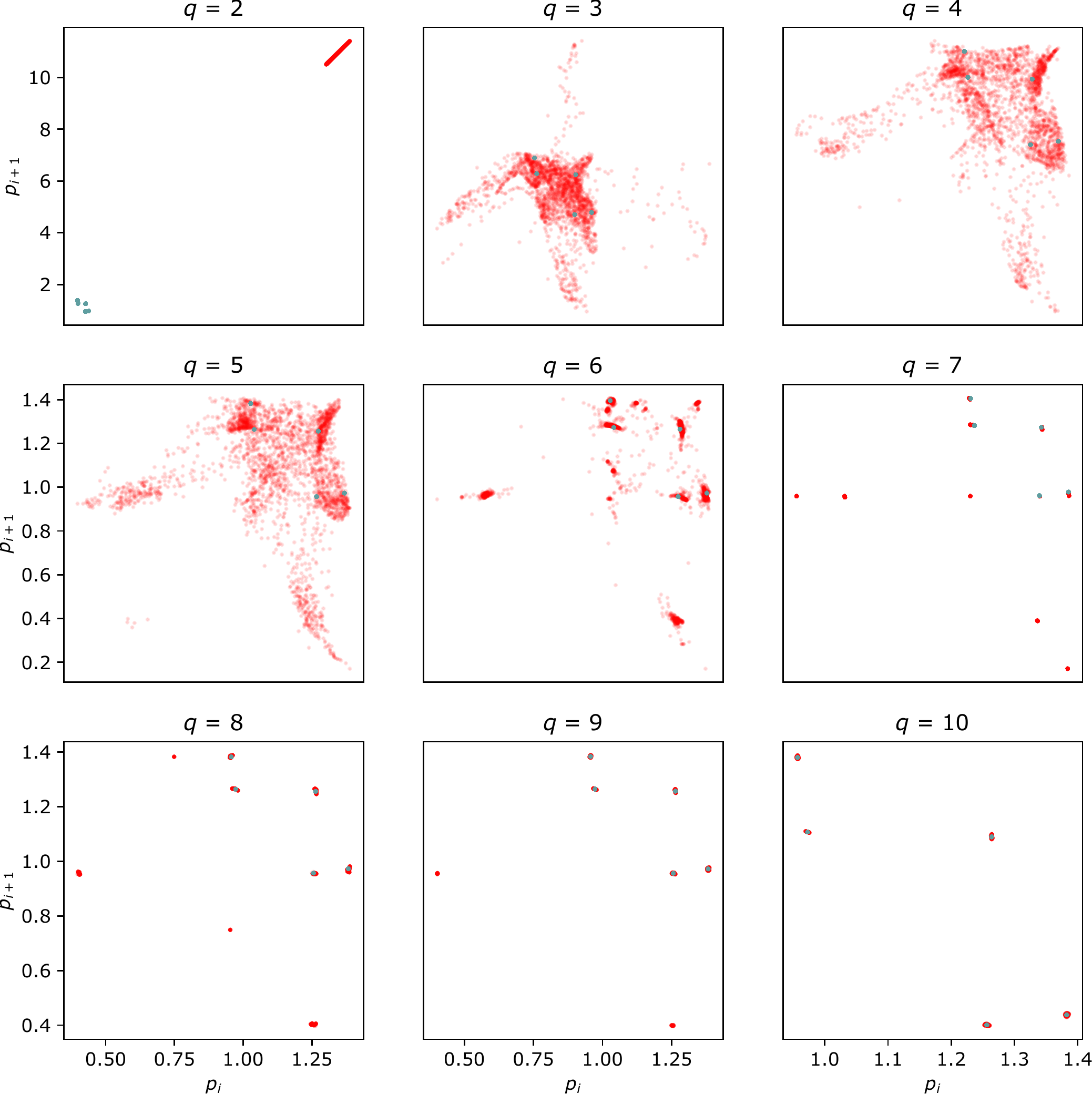}
\caption{Peak-to-peak maps for the ground-truth ($q=17$) Mackey-Glass system with system delay parameters $\tau=4$  and  $\varepsilon=0.1$ (blue dots) compared to approximating system ($q=2,\ldots,10$, red dots). Peaks extracted from solutions over the time interval [0,2.5e3].}
\label{fig:mg_401_p2p_comps}
\end{center}
\end{figure}

Analogous plots for all choices of systems given in Table \ref{tab:numexp} are provided in the supplement, and showcase the significant dependence of dynamics on $q$, at least for small $q$.  That said, in all cases explored here, there is eventually good agreement between the peak-to-peak maps for $q \leq 10$ (often for much smaller $q$) and the ground truth solutions. 

To quantify the apparent convergence/stabilization of approximate solutions across different numbers of interpolating nodes, we constructed delay embeddings of time series samples of each trajectory, after normalizing by the standard deviation of the time series. More precisely, we constructed temporally-ordered point clouds $\{{\bm y_i}\} \subset R^d$ with 
$$\bm y_{i} = [y_i, y_{i+1}, \ldots, y_{i+d-1}],$$
$y_i = y(i-1)/\text{std}(\{y_i\})$, for $i=1, \ldots,$ 1e4, and $d=3$ for the Ikeda systems and the Mackey-Glass systems with $\tau = 2$ or $d=4$ for the Mackey-Glass systems with $\tau = 4$. These embedding dimensions were chosen because they yielded less than 1\% false nearest neighbors using the false nearest neighbor criterion \cite{kennelfnn} with thresholds $R_{\text{tol}}=10$ and $A_{\text{tol}}=2$, and so we expect are sufficiently large to unfold the underlying attractors. Thus, for each simulation a point cloud consisting of approximately 1e4 points in either $\mathbb{R}^3$ or $\mathbb{R}^4$ was constructed, ostensibly capturing a sampling of a diffeomorphic copy of the underlying attractor. 

As a first, coarse comparison of dynamic behavior we estimated the correlation dimension (CD) \cite{corrdim} of the apparent underlying attractor of each system given in Table \ref{tab:numexp}. For details see Section \ref{ssec:quantmethods_cd}. The median and interquartile ranges of the estimated CD are given in Tables \ref{tab:quant_metrics_cubic}, \ref{tab:quant_metrics_mg_2} and \ref{tab:quant_metrics_mg_4}.

We further refined the comparison of attractor structure beyond dimension using persistence homology to quantify topological structures of the attractors.  In particular we computed quantified representations (persistence diagrams) of the connectivity, hole, and void structures of the delay embeddings. For details see Section \ref{ssec:quantmethods_ph}. The median and interquartile ranges of the Wasserstein distances between $H_0$, $H_1$, and $H_2$ diagrams are given in Tables \ref{tab:quant_metrics_cubic}, \ref{tab:quant_metrics_mg_2} and \ref{tab:quant_metrics_mg_4}.

For the Ikeda cubic system with $\tau=1.62$ and $\varepsilon=0$ the quantitative metrics reinforce the qualitative observations made using peak-to-peak maps. For instance, for the Ikeda system with $\varepsilon=0$, the apparent stabilization of the estimated correlation dimension and topological structure occurs at $q=4$ (Figure \ref{fig:cubic1620}). The discrepancies between these structural features for $q<4$ and $q \geq 4$ is pronounced across each metric, except perhaps the Wasserstein distances between $H_2$ diagrams---although we note that the median Wasserstein distances for the $H_2$ diagrams do, in fact, fall outside the baseline interquartile ranges only for $q < 4$. We expect this diminished convergence signal is due to the absence of any high-persistence $H_2$ classes, as can be seen in the $H_2$ diagrams (Figure \ref{fig:cubic1620} (b)).  This indicates the attractors do not exhibit high-persistence void structures, which is consistent with the Lissajou plots of the ground-truth system (Figure \ref{fig:cubic1620}) and the estimated correlation dimensions (Table \ref{tab:quant_metrics_cubic}).

The sudden stabilization of attractor structural features at a small number of interpolating nodes appears to be a general phenomenon across a range of parameters (See Tables \ref{tab:quant_metrics_cubic}, \ref{tab:quant_metrics_mg_2} and \ref{tab:quant_metrics_mg_4}). However, we observe a variety of behaviors for different systems and system delay parameters as $q$ increases. For example, for the Mackey-Glass map with $\tau=4$ and $\varepsilon = 0.1$ we observe notably higher correlation dimensions for $q=3,4,5,6$ than the ground truth system and higher-persistence $H_2$ classes, which is reflected in the Wasserstein distances between the $H_2$ diagrams in this range of $q$ and the ground truth diagrams being much larger than might be expected due to topological noise and finite sampling (Figure \ref{fig:mg401}). The same observation holds for the Ikeda system with $\varepsilon = 0.15$ and $q=3, 4$. Interestingly, for $q=5$, the system appears to collapse to very nearly a periodic orbit. For both of these cases, the true system appears to be 1-dimensional topological circle (Figures \ref{fig:mg_ground_4} and \ref{fig:cubic_ground_truth}), with stabilization occurring around $q=9$ for Mackey-Glass (Table \ref{tab:quant_metrics_mg_4}) and $q=7$ for Ikeda (Table \ref{tab:quant_metrics_cubic}).

\begin{figure}[ht]
\begin{center}
\includegraphics[width=1.0\textwidth]{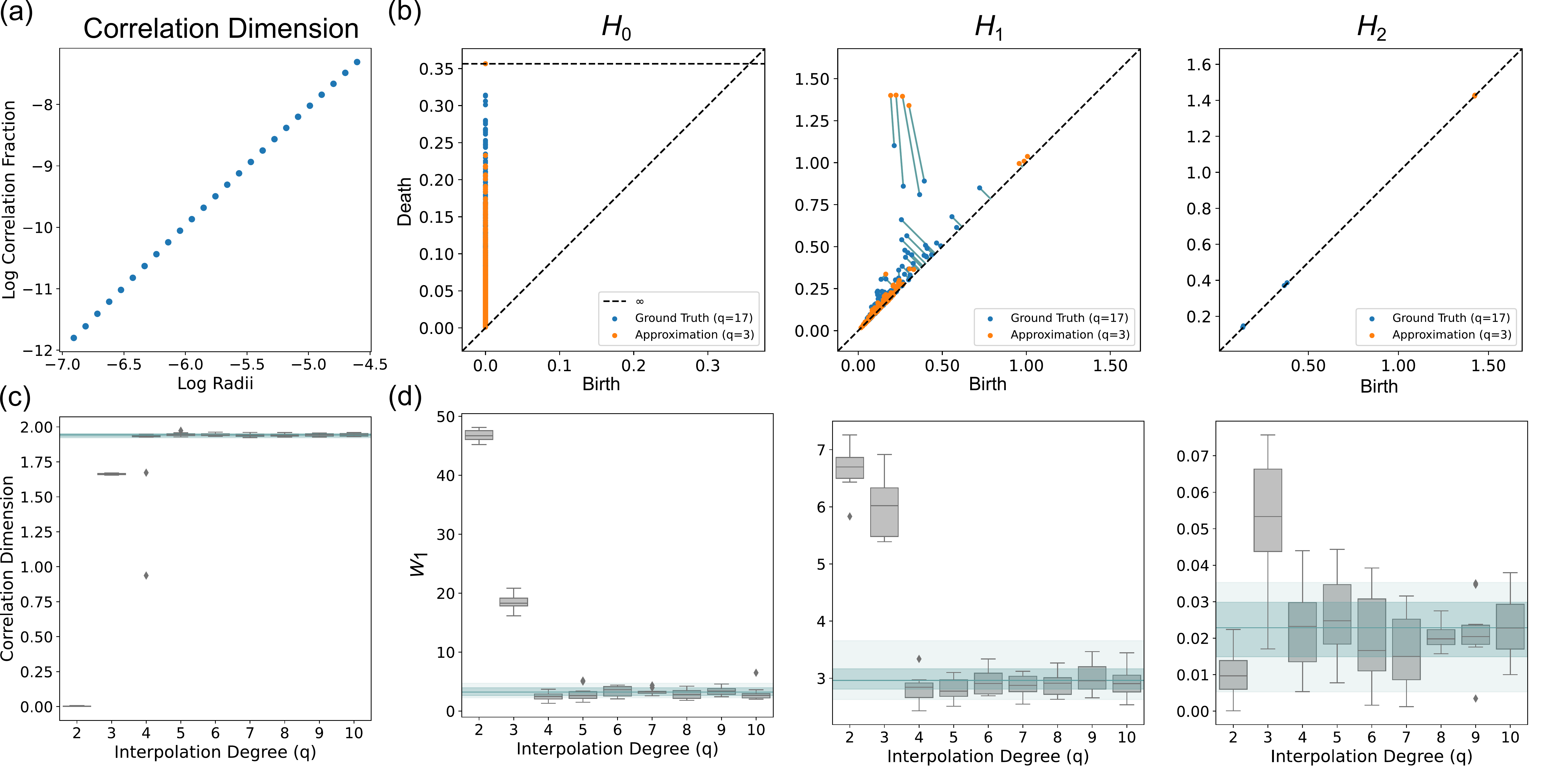}
\caption{Quantitative metrics for the Ikeda cubic system with system delay parameters $\tau=1.62$ and $\varepsilon=0$. (a) Log of correlation fraction ($C_i$) against log of radii parameter ($r_i$) for a single sampling of the ground truth attractor ($q=17$).  Example of data used to estimate the correlation dimension as the slope of the line of best fit. (b) $H_0$, $H_1$, and $H_2$ persistence diagrams for a single sampling of the ground truth attractor ($q=17$) and the corresponding $q=3$ approximation. The persistence pair matching which achieves the $W_1$ distance is shown for the $H_1$ diagrams. (c) Box plots of correlation dimension estimates over multiple samplings of the attractors for different numbers of approximations of interpolating nodes. Also shown is the median (solid horizontal line), interquartile range (dark shaded region), and extent (light shaded region) of the correlation dimension estimates for samplings of the ground truth attractor. (d) Box plots of 1-Wasserstein distances between $H_0$, $H_1$, and $H_2$ diagrams over multiple samplings of the attractors for different numbers of interpolating nodes. Also shown are the median (solid horizontal line), interquartile range (dark shaded region), and maximum extent (light shaded region) of the 1-Wasserstein distances between different samplings of the ground truth attractor.}
\label{fig:cubic1620}
\end{center}
\end{figure}

\begin{figure}[ht]
\begin{center}
\includegraphics[width=1.0\textwidth]{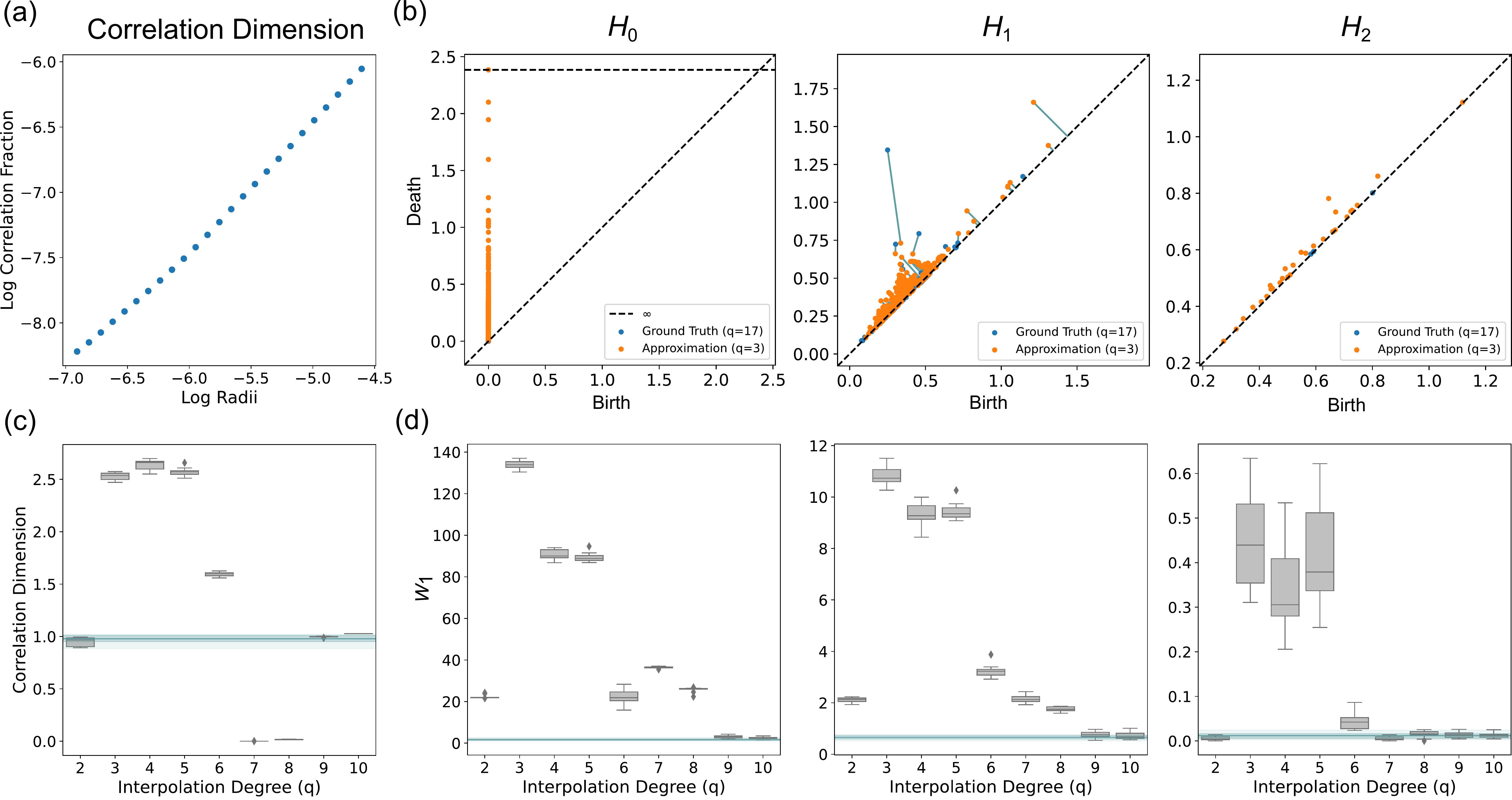}
\caption{Quantitative metrics for the Mackey-Glass system with system delay parameters $\tau=4$ and $\varepsilon=0.1$. (a) Log of correlation fraction ($C_i$) against log of radii parameter ($r_i$) for a single sampling of the ground truth attractor ($q=17$).  Example of data used to estimate the correlation dimension as the slope of the line of best fit. (b) $H_0$, $H_1$, and $H_2$ persistence diagrams for a single sampling of the ground truth attractor ($q=17$) and the corresponding $q=3$ approximation. The persistence pair matching which achieves the $W_1$ distance is shown for the $H_1$ diagrams. (c) Box plots of correlation dimension estimates over multiple samplings of the attractors for different numbers of approximations of interpolating nodes. Also shown is the median (solid horizontal line), interquartile range (dark shaded region), and extent (light shaded region) of the correlation dimension estimates for samplings of the ground truth attractor. (d) Box plots of 1-Wasserstein distances between $H_0$, $H_1$, and $H_2$ diagrams over multiple samplings of the attractors for different numbers of interpolating nodes. Also shown are the median (solid horizontal line), interquartile range (dark shaded region), and maximum extent (light shaded region) of the 1-Wasserstein distances between different samplings of the ground truth attractor.}
\label{fig:mg401}
\end{center}
\end{figure}

\begin{landscape}
\thispagestyle{empty}
\begin{table}
\centering
\small
\begin{tabular}{llllllllllll}
\toprule
    &       &             2 &            3 &            4 &            5 &            6 &           7 &           8 &           9 &          10 &    baseline \\
eps & metric &               &              &              &              &              &             &             &             &             &             \\
\midrule
\multirow{4}{*}{0} & Dim. &    0.00 (0.0) &   1.66 (0.0) &   1.93 (0.0) &   1.94 (0.0) &   1.94 (0.0) &  1.94 (0.0) &  1.94 (0.0) &  1.94 (0.0) &  1.95 (0.0) &  1.94 (0.0) \\
    & $H_0$ &   46.71 (1.5) &  18.28 (1.3) &   2.48 (0.9) &   2.61 (1.1) &   3.63 (1.6) &  3.23 (0.4) &  2.78 (1.3) &  3.38 (1.0) &  2.65 (0.8) &  3.22 (1.3) \\
    & $H_1$ &    6.70 (0.4) &   6.02 (0.9) &   2.84 (0.3) &   2.78 (0.3) &   2.91 (0.4) &  2.88 (0.3) &  2.91 (0.3) &  2.96 (0.4) &  2.91 (0.3) &  2.96 (0.4) \\
    & $H_2$ &    0.01 (0.0) &   0.05 (0.0) &   0.02 (0.0) &   0.02 (0.0) &   0.02 (0.0) &  0.02 (0.0) &  0.02 (0.0) &  0.02 (0.0) &  0.02 (0.0) &  0.02 (0.0) \\
\cline{1-12}
\multirow{4}{*}{0.05} & Dim. &    0.02 (0.0) &   1.88 (0.0) &   1.94 (0.0) &   1.92 (0.0) &   1.94 (0.0) &  1.94 (0.0) &  1.94 (0.0) &  1.94 (0.0) &  1.94 (0.0) &  1.94 (0.0) \\
    & $H_0$ &  41.64 (15.8) &   8.08 (2.4) &   2.87 (0.8) &   2.95 (0.9) &   2.63 (1.3) &  2.59 (1.0) &  2.69 (0.9) &  3.25 (1.1) &  2.75 (0.9) &  3.43 (1.5) \\
    & $H_1$ &    5.37 (0.9) &   3.04 (0.2) &   2.49 (0.1) &   2.53 (0.2) &   2.40 (0.1) &  2.43 (0.2) &  2.54 (0.2) &  2.42 (0.2) &  2.57 (0.2) &  2.51 (0.2) \\
    & $H_2$ &    0.01 (0.0) &   0.01 (0.0) &   0.02 (0.0) &   0.02 (0.0) &   0.02 (0.0) &  0.02 (0.0) &  0.03 (0.0) &  0.02 (0.0) &  0.02 (0.0) &  0.01 (0.0) \\
\cline{1-12}
\multirow{4}{*}{0.15} & Dim. &  - &   1.73 (0.0) &   1.21 (0.0) &   0.00 (0.0) &   0.98 (0.0) &  1.01 (0.0) &  1.02 (0.0) &  1.02 (0.0) &  1.01 (0.0) &  1.01 (0.0) \\
    & $H_0$ &   -  &  25.93 (0.6) &   2.96 (0.3) &  20.87 (0.3) &  10.17 (0.4) &  0.84 (0.1) &  0.87 (0.1) &  0.86 (0.2) &  0.62 (0.3) &  0.68 (0.2) \\
    & $H_1$ &    -  &   2.45 (0.1) &   0.22 (0.1) &   1.04 (0.0) &   0.41 (0.0) &  0.04 (0.0) &  0.04 (0.0) &  0.05 (0.0) &  0.03 (0.0) &  0.03 (0.0) \\
    & $H_2$ &   -  &   0.01 (0.0) &   0.00 (0.0) &   0.00 (0.0) &   0.00 (0.0) &  0.00 (0.0) &  0.00 (0.0) &  0.00 (0.0) &  0.00 (0.0) &  0.00 (0.0) \\
\cline{1-12}
\multirow{4}{*}{0.2} & Dim. &    1.00 (0.0) &   - &   - &   1.82 (0.0) &   1.74 (0.0) &  1.77 (0.0) &  1.79 (0.0) &  1.79 (0.0) &  1.79 (0.0) &  1.79 (0.0) \\
    & $H_0$ &   39.22 (1.7) &  - &  - &   3.55 (2.0) &   3.38 (0.8) &  3.32 (1.9) &  3.35 (1.6) &  3.13 (0.8) &  3.18 (1.6) &  4.25 (1.8) \\
    & $H_1$ &    4.54 (0.7) &   - &   - &   2.80 (0.4) &   2.66 (0.6) &  2.71 (0.2) &  2.71 (0.5) &  2.54 (0.5) &  2.79 (0.6) &  2.80 (0.2) \\
    & $H_2$ &    0.02 (0.0) &   - &   - &   0.03 (0.0) &   0.03 (0.0) &  0.03 (0.0) &  0.03 (0.0) &  0.03 (0.0) &  0.04 (0.0) &  0.04 (0.0) \\
\bottomrule
\end{tabular}
\caption{\textbf{Quantitative Comparisons of Attractors, Cubic Ikeda}, \bm{$\tau=1.62$}: For each choice of  state-dependent delay parameter $\varepsilon$, the top rows provide median and interquartile ranges of estimated correlation dimensions (Dim.) across different numbers of interpolating nodes ($q=2,\ldots 10$) and estimates for $q=17$ in the baseline column. Rows 2,3, and 4 give median and interquartile ranges of the 1-Wasserstein distance between $H_0$, $H_1$, and $H_2$ persistence diagrams between the approximate systems ($q<17$) and the ground truth system ($q=17$), with baseline median and interquratile ranges measuring 1-Wasserstein distances between diaagrams corresponding to different random samplings of the ground truth ($q=17$) attractor. Missing values due to numerical instabilities for these choices of $q$.}
\label{tab:quant_metrics_cubic}
\end{table}
\end{landscape}

\begin{landscape}
\thispagestyle{empty}
\begin{table}
\centering
\small
\begin{tabular}{llllllllllll}
\toprule
     &       &            2 &            3 &           4 &           5 &           6 &           7 &           8 &           9 &          10 &    baseline \\
eps & metric &              &              &             &             &             &             &             &             &             &             \\
\midrule
\multirow{4}{*}{0} & Dim. &   0.02 (0.0) &   1.81 (0.0) &  1.94 (0.0) &  2.09 (0.0) &  2.09 (0.0) &  2.09 (0.0) &  2.08 (0.0) &  2.09 (0.0) &  2.08 (0.0) &  2.10 (0.0) \\
     & $H_0$ &  73.57 (1.1) &  18.59 (1.2) &  3.82 (1.6) &  2.58 (0.6) &  3.14 (0.7) &  3.76 (1.3) &  2.99 (0.8) &  2.77 (1.3) &  3.15 (0.9) &  3.24 (0.9) \\
     & $H_1$ &   9.65 (0.3) &   4.42 (0.5) &  3.15 (0.2) &  3.00 (0.3) &  2.95 (0.3) &  3.00 (0.3) &  3.23 (0.6) &  3.25 (0.6) &  2.97 (0.4) &  3.10 (0.6) \\
     & $H_2$ &   0.18 (0.1) &   0.17 (0.1) &  0.21 (0.1) &  0.25 (0.1) &  0.21 (0.0) &  0.22 (0.0) &  0.23 (0.1) &  0.21 (0.1) &  0.20 (0.1) &  0.21 (0.0) \\
\cline{1-12}
\multirow{4}{*}{0.05} & Dim. &  -0.00 (0.0) &   1.52 (0.0) &  1.96 (0.0) &  2.08 (0.0) &  2.11 (0.0) &  2.12 (0.0) &  2.11 (0.0) &  2.10 (0.0) &  2.10 (0.0) &  2.11 (0.0) \\
     & $H_0$ &  64.46 (0.8) &  33.29 (1.8) &  4.22 (1.1) &  2.38 (0.8) &  2.48 (1.1) &  2.44 (0.8) &  2.27 (0.6) &  1.99 (1.0) &  2.48 (0.9) &  2.14 (0.5) \\
     & $H_1$ &   8.56 (0.3) &   6.65 (0.5) &  2.73 (0.3) &  2.63 (0.4) &  2.71 (0.2) &  2.75 (0.3) &  2.67 (0.4) &  2.78 (0.3) &  2.63 (0.3) &  2.83 (0.4) \\
     & $H_2$ &   0.10 (0.1) &   0.07 (0.1) &  0.11 (0.0) &  0.09 (0.0) &  0.11 (0.0) &  0.09 (0.0) &  0.10 (0.0) &  0.10 (0.0) &  0.09 (0.0) &  0.11 (0.1) \\
\cline{1-12}
\multirow{4}{*}{0.15} & Dim. &   1.70 (0.0) &   1.54 (0.0) &  1.99 (0.0) &  2.06 (0.0) &  2.04 (0.0) &  2.04 (0.0) &  2.02 (0.0) &  2.03 (0.0) &  2.03 (0.0) &  2.02 (0.0) \\
     & $H_0$ &  12.58 (3.4) &  12.59 (2.1) &  2.41 (0.8) &  2.25 (0.9) &  2.20 (0.4) &  2.05 (0.5) &  2.12 (0.7) &  2.02 (0.5) &  2.70 (0.5) &  2.49 (0.6) \\
     & $H_1$ &   3.75 (0.4) &   3.44 (0.6) &  2.42 (0.1) &  2.32 (0.2) &  2.30 (0.3) &  2.33 (0.2) &  2.35 (0.4) &  2.31 (0.1) &  2.34 (0.3) &  2.34 (0.2) \\
     & $H_2$ &   0.03 (0.0) &   0.03 (0.0) &  0.03 (0.0) &  0.04 (0.0) &  0.03 (0.0) &  0.04 (0.0) &  0.04 (0.0) &  0.03 (0.0) &  0.03 (0.0) &  0.03 (0.0) \\
\bottomrule
\end{tabular}
\caption{\textbf{Quantitative Comparisons of Attractors, Mackey-Glass}, \bm{$\tau=2$}: For each choice of  state-dependent delay parameter $\varepsilon$, the top rows provide median and interquartile ranges of estimated correlation dimensions (Dim.) across different numbers of interpolating nodes ($q=2,\ldots 10$) and estimates for $q=17$ in the baseline column. Rows 2,3, and 4 give median and interquartile ranges of the 1-Wasserstein distance between $H_0$, $H_1$, and $H_2$ persistence diagrams between the approximate systems ($q<17$) and the ground truth system ($q=17$), with baseline median and interquratile ranges measuring 1-Wasserstein distances between diaagrams corresponding to different random samplings of the ground truth ($q=17$) attractor.}
\label{tab:quant_metrics_mg_2}
\end{table}
\end{landscape}

\begin{landscape}
\thispagestyle{empty}
\begin{table}
\centering
\small
\begin{tabular}{llllllllllll}
\toprule
    &       &             2 &             3 &             4 &            5 &            6 &            7 &            8 &             9 &           10 &    baseline \\
eps & metric &               &               &               &              &              &              &              &               &              &             \\
\midrule
\multirow{4}{*}{0} & Dim. &   -0.00 (0.0) &    1.02 (0.0) &    3.17 (0.1) &   3.59 (0.3) &   3.64 (0.5) &   3.36 (0.3) &   3.56 (0.4) &    3.45 (0.2) &   3.55 (0.6) &  3.42 (0.4) \\
    & $H_0$ &  173.48 (2.5) &  179.31 (7.4) &   15.29 (6.8) &   7.03 (4.8) &   5.86 (3.0) &   5.38 (1.1) &   5.20 (2.7) &    6.43 (1.7) &   6.91 (1.7) &  6.26 (3.8) \\
    & $H_1$ &   19.60 (1.1) &   20.63 (1.2) &    6.83 (1.7) &   5.82 (0.9) &   5.83 (0.6) &   5.76 (0.6) &   5.80 (0.8) &    5.79 (0.6) &   6.25 (0.8) &  5.90 (0.5) \\
    & $H_2$ &    1.62 (0.2) &    1.63 (0.2) &    1.48 (0.2) &   1.36 (0.1) &   1.36 (0.2) &   1.26 (0.3) &   1.38 (0.3) &    1.30 (0.2) &   1.28 (0.1) &  1.28 (0.2) \\
\cline{1-12}
\multirow{4}{*}{0.05} & Dim. &    1.01 (0.0) &    2.45 (0.1) &    2.46 (0.1) &   2.95 (0.1) &   2.51 (0.1) &   2.18 (0.0) &   2.05 (0.0) &    1.72 (0.1) &   1.90 (0.3) &  2.12 (0.0) \\
    & $H_0$ &   59.26 (3.5) &  82.25 (12.3) &  90.36 (10.5) &  94.19 (8.5) &  40.90 (6.2) &  16.58 (5.3) &   7.03 (3.0) &  23.98 (16.7) &  9.03 (10.1) &  8.46 (8.0) \\
    & $H_1$ &    5.52 (0.5) &   11.10 (1.9) &   12.03 (1.9) &  12.15 (1.2) &   5.97 (0.6) &   3.44 (0.4) &   2.67 (0.8) &    3.91 (1.2) &   2.86 (1.2) &  3.07 (1.2) \\
    & $H_2$ &    0.03 (0.0) &    0.63 (0.1) &    0.73 (0.2) &   0.73 (0.1) &   0.21 (0.1) &   0.08 (0.0) &   0.07 (0.0) &    0.04 (0.0) &   0.06 (0.0) &  0.07 (0.0) \\
\cline{1-12}
\multirow{4}{*}{0.1} & Dim. &    0.97 (0.1) &    2.54 (0.1) &    2.66 (0.1) &   2.57 (0.0) &   1.60 (0.0) &   0.00 (0.0) &   0.02 (0.0) &    1.00 (0.0) &   1.03 (0.0) &  0.98 (0.1) \\
    & $H_0$ &   21.92 (0.2) &  133.94 (2.7) &   90.04 (3.9) &  88.88 (2.3) &  21.84 (4.2) &  36.35 (0.4) &  26.14 (0.3) &    3.15 (0.9) &   2.45 (0.5) &  1.71 (0.8) \\
    & $H_1$ &    2.14 (0.1) &   10.74 (0.5) &    9.28 (0.5) &   9.35 (0.4) &   3.23 (0.2) &   2.13 (0.2) &   1.75 (0.2) &    0.76 (0.2) &   0.68 (0.2) &  0.65 (0.1) \\
    & $H_2$ &    0.00 (0.0) &    0.44 (0.2) &    0.31 (0.1) &   0.38 (0.2) &   0.04 (0.0) &   0.00 (0.0) &   0.02 (0.0) &    0.01 (0.0) &   0.01 (0.0) &  0.01 (0.0) \\
\bottomrule
\end{tabular}
\caption{\textbf{Quantitative Comparisons of Attractors, Mackey-Glass}, \bm{$\tau=4$}: For each choice of  state-dependent delay parameter $\varepsilon$, the top rows provide median and interquartile ranges of estimated correlation dimensions (Dim.) across different numbers of interpolating nodes ($q=2,\ldots 10$) and estimates for $q=17$ in the baseline column. Rows 2,3, and 4 give median and interquartile ranges of the 1-Wasserstein distance between $H_0$, $H_1$, and $H_2$ persistence diagrams between the approximate systems ($q<17$) and the ground truth system ($q=17$), with baseline median and interquratile ranges measuring 1-Wasserstein distances between diaagrams corresponding to different random samplings of the ground truth ($q=17$) attractor.}
\label{tab:quant_metrics_mg_4}
\end{table}
\end{landscape}

\section{Conclusions}

In this work we implemented and profiled an iterative and
perturbative numerical scheme for computing 
orbits of state dependent perturbations of delay maps.
After showing local convergence of the algorithm, we also explored its more qualitative behaviors. Our empirical investigations suggest that the estimated correlation dimensions, the peak-to-peak, and Lissajou plots 
of the analyzed systems support low-dimensional chaotic dynamics across a range of state-dependent delays, and 
remind us of those infinite-dimensional systems known to support inertial manifolds. We observe also windows of the 
state-dependent delay parameter on which the apparent chaotic dynamics of the state-independent delay system 
appear to collapse to periodic orbits.

We observe apparent convergence of qualitative and quantitative features of the dynamics for relatively small numbers of interpolating nodes (\(q\leq 10\)) 
across systems and system parameters. Furthermore, we observe a rich set of behaviors for small 
values of \( q \), 
in contrast to the ground truth systems. For instance, for sufficiently large \(q\) 
the system may be 
periodic, behaving qualitatively like the true flow the map represents, while for smaller \(q\) 
the map can appear to exhibit chaotic dynamics. 

 We believe that the above observations hold for a large class of systems 
similar to those under consideration in this paper; that is, when the 
right hand side of the DDE depends in a linear manner with respect to the 
instantaneous state variable \( x(t)\)
(even if the dependence in the delayed variable 
is more complex).
We also believe that those observations extend to more 
general systems where the solution of Eq. (6), 
being a fixed point of a specific Picard operator, can
be computed thanks to a Newton-like method. 
This question will be explored in a forthcoming project.

\section{Methods}

\subsection{Quantitative Attractor Metrics}
\label{sec:quantmethods}

\subsubsection{Correlation Dimension}
\label{ssec:quantmethods_cd}
In \cite{corrdim} it was argued that the correlation integral,
\[
C(l) := \lim_{N \rightarrow \infty} \# \{ (m, n) \; | \; \|{\bm y}_m-{\bm y}_n \| < l\} / N^2,
\]
grows according to a power law
\[
C(l) \sim l^{\nu},
\]
when \({\bm y}_m\) 
are points on an attractor associated with a nonlinear dissipative dynamical system sampled at fixed time increments.

It was proposed 
that the CD, \(\nu\), 
serves as a measure of 
dynamic complexity 
and is related to other measures such as the Hausdorff \cite{Hausdorff1918} and information dimensions \cite{infdim}, while being more readily computable using standard linear regression, in the log-log plane, of estimates of $C(l)$ for small $l$.

For each choice of system parameters in this study and for each number of interpolating nodes, $q$, fifty simulations were performed with random initial conditions. From these, ten sets of five embedded point clouds were concatenated to very densely sample each system's attractor. The CD of each resulting point cloud was then estimated---using standard linear regression---as the slope of the line of best fit through twenty-five data points, $(\log(C_i), \log(r_i))$, with $r_i$ being logarithmically spaced in the interval [1e-3, 1e-2] and 
$$C_i = 2\# \{ (m, n) \; | \; 1\leq m < n \leq N, \|{\bm y}_m-{\bm y}_n \| < r_i\} / N(N-1)$$
being the fraction of Euclidean distances between distinct points in the cloud that are less than $r_i$. We refer to $C_i$ as the correlation fraction corresponding to the radius parameter $r_i$.

\subsubsection{Persistent Homology}
\label{ssec:quantmethods_ph}
In addition to CD, we computed three topological descriptors of the reconstructed attractors using persistent homology (PH). Recently PH has found widespread use across experimental sciences and in machine learning applications for its ability to extract and represent quantified shape-based features. For completeness we provide a brief description of the particular use of PH in this study. For a more complete treatment of the subject and its applications see one of the recent review articles \cite{Otter2017, adams2021review, chazal2021roadmap}. 

PH may be regarded as a transformation mapping a finite point-cloud ${\bm Y} := \{{\bm y}_m\} \subset \mathbb{R}^d$ to a collection of so-called persistence diagrams (PD) of the form
$$\text{dgm}_k({\bm Y}) = \{(b_i, d_i) \; | \; b_i < d_i \in \mathbb{R} \} \cup \prod_{n=1}^{\infty} \Delta,$$ 
where $\Delta := \{(x,x) \; | \; x \in \mathbb{R}\}$ is the diagonal in the plane, $k\geq 0$. Thus, each PD consists of a finite multiset of planar points above the diagonal with countably infinitely many copies of the diagonal. We refer to elements of a PD as persistence pairs. 

The inclusion of copies of the diagonal allows us to define the $p$-Wasserstein distance ($p \geq 1$) between PDs: 
$$
W_p(\text{dgm}_{k}({\bm Y}), \text{dgm}_{k}({\bm X})) := \inf_{\phi}\left(\sum_{{\bm a} \in \text{dgm}_{k}({\bm Y})} \|{\bm a} - \phi({\bm a})\|^p\right)^{1/p},
$$ 
where the infimum is taken over all bijections $\phi: \text{dgm}_{k}({\bm Y}) \rightarrow \text{dgm}_{k}({\bm Y})$. Intuitively, $W_p$ measures the minimal sum of distances between matchings of persistence pairs in two diagrams among all possible pairings of those points, allowing off-diagonal pairs in one diagram to pair to diagonal points in the other. Although the sum is technically over a countably-infinite multiset, only off-diagonal pairs will contribute to the sum for a bijection which realizes the infimum, which must exist for diagrams with finitely many off-diagonal pairs.

A persistence diagram can be constructed from a finite simplicial complex $K$ by defining a filtration on it,
$\emptyset = K_0 \subseteq K_1 \subseteq \ldots \subseteq K_n = K,$ 
thought of as a nested family of simplicial complexes. Computing dimension-$k$ simplicial homology over a field $\mathbb{K}$ on each complex produces a persistence module
$$
\textbf{H}_k(K) := (0 = H_k(K_0) \xrightarrow[]{f_k^0} H_k(K_1) \xrightarrow[]{f_k^1} \cdots \xrightarrow[]{f_k^{n-1}} H_k(K_n) = H_k(K)),
$$
comprised of vector spaces over $\mathbb{K}$ connected by linear maps $$f_{k}^{i,j} := f_k^{j-1}\circ\cdots\circ f_k^i: H_k(K_i) \rightarrow H_k(K_j),$$
induced by inclusion at the level of the simplicial complexes for each $i \leq j$.
$\textbf{H}_k(K)$ can be decomposed into a direct sum, 
$$\bigoplus_{[a, b) \in \mathcal{D}_k} I_{[a, b)},$$
of so-called interval (persistence) modules of the form
$$
I_{[a, b)} = (V_1 \rightarrow \cdots \rightarrow V_n),
$$
where 
$$
V_i = \begin{cases}\mathbb{K},&  i \in [a,b)\\ 0 , & \text{otherwise} \end{cases}
$$
that are connected by maps $f^{i,j}: V_i \rightarrow V_j$ where 
$f^{i,j}$ is the identity map if $a \leq i \leq j < b$ and identically 0 otherwise, and the collection of intervals $\mathcal{D}_k$ are uniquely determined by the filtration. 

Intuitively each interval $[a, b)$ in the interval-module representation encodes a $k$-dimensional homology class that is born (first appears) when moving from some $\text{K}_{i-1}$ to $\text{K}_{i}$ and which dies (becomes homologous to a class which was born earlier) when moving from $\text{K}_{j-1}$ to $\text{K}_{j}$. We refer to the value $b-a$ as the persistence or lifetime of the class corresponding to the interval $[a,b)$. Note, some classed may never die, giving rise to an unbounded interval $[a,\infty)$, but we will omit these from diagrams and calculation of $W_p$ between diagrams. The algebraic information in the interval modules is encoded into a persistence diagram $\text{dgm}_k(K) = \{(a,b) \; | \; [a,b) \in D_k, a,b \in \mathbb{R}\}$ where the birth and death filtration values of each finite-persistence $k$-dimensional homological feature appears and later dies as a result of the filtration.

Here we adopt the commonly used Vietoris-Rips (VR) filtration to construct a filtered simplicial complex from a finite point cloud ${\bm Y} = \{{\bm y}_m\}$, sampled from an attractor. In particular, let $\partial_{i}$ be the $i$-th smallest distance between any two points in ${\bm Y}$ and define
$$\text{VR}_i({\bm Y}) := \{\sigma \subset {\bm Y} \; | \; \|{\bm y}_i-{\bm y}_j \| \leq \partial_i, \text{ for all } {\bm y}_i, {\bm y}_j \in \sigma \}$$
to be the simplicial complex consisting of all subsets whose diameter is less than $\partial_i$. 

We regard this filtration as starting with $\text{VR}_1 = \{ \{{\bm y}_m\} | {\bm y}_m \in {\bm Y}\}$ (at $\partial_1 = 0$), consisting only of 0-dimensional simplices corresponding to the points in the cloud. Edges are added between pairs of points in order from nearest to farthest, with each higher dimensional simplex appearing at the filtration value corresponding to the diameter of the points corresponding to its vertices. In this way, we construct a nested family of complexes meant to capture the intrinsic multiscale geometry of the point cloud. 

Crucial to the use of persistent homology in applications are the numerous stability results which establish---under a variety of assumptions about the data and the metrics placed on data and the diagrams---the (Lipschitz) continuity of these transformations sending data to persistence diagrams \cite{Chazal2012PersistenceSF} \cite{Cohen-SteinerEHM10}. Recently, general stability results for $W_p$, and partial results for $W_p$ when applied to the special case of the VR filtration used here were established \cite{skraba2021wasserstein}.

For each choice of system parameters in this study and for each number of interpolating nodes, 1000 points were randomly subselected from each of ten attractor point clouds generated from ten different initial conditions. $H_0$, $H_1$, and $H_2$ persistence diagrams were computed from the the VR filtered simplicial complexes built on each set of 1000 points using the Ripser algorithm \cite{Ripser, Tralie2018}. The 1-Wasserstein distances between the diagrams derived from each approximate system ($q=2,\ldots, 10$) and the ground truth system ($q=17$) were then computed. To serve as a baseline for comparison, and to account for topological noise due to finite sampling of the attractors, the 1-Wasserstein distances between ten pairs of ground truth diagrams were also computed.

$\ $



\end{document}